\documentclass[preprint,article,accept,moreauthors,pdftex]{Definitions/mdpi}
\firstpage{1} 
\pubvolume{xx}
\issuenum{1}
\articlenumber{1}
\pubyear{2020}
\copyrightyear{2020}
\history{Received: dd.mm.yyyy; Accepted: dd.mm.yyyy ; Published: dd.mm.yyyy}
\usepackage{algorithm2e}
\usepackage{amsmath}
\usepackage{amsfonts}
\usepackage{amsthm}
\usepackage{authblk}
\usepackage[ngerman, UKenglish]{babel}
\usepackage{bbm}
\usepackage{bigstrut}
\usepackage{caption}
\usepackage{cases}
\usepackage{comment}
\usepackage{float}
\usepackage[acronym]{glossaries}
\glsdisablehyper
\usepackage{graphicx, subfigure}
\usepackage[utf8]{inputenc}
\usepackage{mathtools}
\usepackage{multirow}
\usepackage{pgfplots}
\usepackage{placeins}
\usepackage{siunitx}
\usepackage{tabularx}
\usepackage{tikz}

\makeglossaries
\DeclareMathAlphabet\mathbfcal{OMS}{cmsy}{b}{n}
\newacronym{em}{EM}{electromagnetic}
\newacronym{emf}{EMF}{electromagnetic field}
\newacronym{fe}{FE}{finite element}
\newacronym{fem}{FEM}{finite element method}
\newacronym{gpc}{gPC}{generalized polynomial chaos}
\newacronym{pec}{PEC}{perfect electric conductor}
\newacronym{rms}{RMS}{root-mean-square}
\newacronym{te}{TE}{transverse electric}
\newacronym{tm}{TM}{transverse magnetic}
\newacronym{tem}{TEM}{transverse electric and magnetic}
\newacronym{lar}{LAR}{least angle regression}
\newacronym{lhs}{LHS}{latin hypercube sampling}
\newacronym{ls}{LS}{least squares}
\newacronym{mlmc}{MLMC}{multilevel Monte Carlo}
\newacronym{mc}{MC}{Monte Carlo}
\newacronym{pce}{PCE}{polynomial chaos expansion}
\newacronym{pde}{PDE}{partial differential equation}
\newacronym{pdf}{PDF}{probability density function}
\newacronym[longplural={quantities of interest}]{qoi}{QoI}{quantity of interest}
\newacronym{rv}{RV}{random variable}
\newacronym{td}{TD}{total degree}
\newacronym{tp}{TP}{tensor product}
\newacronym{uq}{UQ}{uncertainty quantification}

\DeclareMathOperator*{\argmax}{argmax}
\Title{Approximation and Uncertainty Quantification of Systems with Arbitrary Parameter Distributions using Weighted Leja Interpolation}

\Author{Dimitrios Loukrezis $^{1,2,*,\dagger}$\orcidA and Herbert De Gersem $^{1,2,\dagger}$\orcidA}%Please carefully check the accuracy of names and affiliations.

\AuthorNames{Dimitrios Loukrezis and Herbert De Gersem}

\address{%
$^{1}$ \quad Institute for Accelerator Science and Electromagnetic Fields (TEMF), Technische Universit\"at Darmstadt, 64289 Darmstadt, Germany; degersem@temf.tu-darmstadt.de\\
$^{2}$ \quad Centre for Computational Engineering, Technische Universit\"at Darmstadt, 64293 Darmstadt, Germany}

\corres{Correspondence: loukrezis@temf.tu-darmstadt.de}

\firstnote{Current address: Schlossgartenstra{\ss}e 8, 64289 Darmstadt, Germany} 
\abstract{Approximation and uncertainty quantification methods based on Lagrange interpolation are typically abandoned in cases where the probability distributions of one or more  {system} parameters are not normal, uniform, or closely related  {distributions}, due to the computational issues that arise when one wishes to define interpolation nodes for general distributions. 
This paper examines the use of the recently introduced weighted Leja nodes for that purpose. 
Weighted Leja interpolation rules are presented, along with a dimension-adaptive sparse interpolation algorithm, to be employed in the case of high-dimensional input uncertainty. 
The performance and reliability of the suggested approach is verified by four numerical experiments, where the respective models feature extreme value and truncated normal parameter distributions. 
 {Furthermore}, the suggested approach is compared   {with} a well-established polynomial chaos method and found to be either comparable or superior in terms of approximation and statistics estimation accuracy.}


\begin{document}
%%%%%%%%%%%%%%%%%%%%%%%%%%%%%%%%%%%%%%%%%%

%%%%%%%%%%%%%%%%%%%%%%%%%%%%%%%%%%%%%%%%%%
\section{Introduction}
\label{sec:intro}
Studying physical systems often involves assessing the influence of random parameters upon a system's behavior and outputs.
Such \gls{uq} studies are becoming increasingly popular in a variety of technical domains. 
Of particular importance to industrial scientists and engineers are {so-called} non-intrusive \gls{uq} methods, where the term ``non-intrusive'' signifies that the often complex and typical proprietary simulation software applied to model physical systems are incorporated into \gls{uq} studies without any modifications.

Sampling methods such as (quasi-) \gls{mc} \cite{caflisch1998} or \gls{lhs} \cite{loh1996} are the most common approaches applied for non-intrusive \gls{uq}, as they are easy to implement, straightforwardly parallelizable, and can be applied to problems with a large number of parameters under relatively mild assumptions regarding the regularity of the underlying mathematical model.
However, their slow convergence rates render them unattractive or even intractable in cases where high accuracy is demanded and single simulations are highly time-consuming.

Computationally efficient alternatives are spectral \gls{uq} methods \cite{lemaitre2010, xiu2010}, which approximate the functional dependence of a system's \gls{qoi} on its input parameters. 
The most popular black-box methods of this category employ global polynomial approximations by aid of either Lagrange interpolation schemes \cite{babuska2010, barthelmann2000} or \gls{gpc} \cite{xiu2002}, where the latter approach is typically based on \gls{ls} regression \cite{blatman2010, blatman2011, migliorati2013a} or pseudo-spectral projection methods \cite{conrad2013, constantine2012, lemaitre2001}, or, less commonly, on interpolation \cite{loukrezis2019b}.
Assuming a smooth input-output dependence, spectral \gls{uq} methods provide fast convergence rates, in some cases even of exponential order.
However, their performance decreases rapidly for an increasing number of input parameters due to the curse of dimensionality \cite{bellman1957}.
To overcome this bottleneck, state-of-the-art approaches rely on algorithms for the adaptive construction of sparse approximations \cite{blatman2010, blatman2011, chkifa2014, migliorati2013a, nobile2008a, stoyanov2016}.

Adaptive sparse approximation algorithms for Lagrange interpolation methods typically employ nested sequences of univariate interpolation nodes. 
{While not strictly necessary \cite{ernst2018}, nested node sequences are helpful for the efficient construction of sparse grid interpolation algorithms.
Nested node sequences are readily available for the  well studied cases of uniformly or normally distributed \glspl{rv} \cite{clenshaw1960, genz1996}, but not for more general probability measures.
In principle, one could derive nested interpolation (or quadrature) rules tailored to an arbitrary \gls{pdf} \cite{burkardt2014, sommariva2013}, however, this is a rather cumbersome task for which dedicated analyses are necessary each and every time a new \gls{pdf} is considered.}

Diversely, the main requirement in \gls{gpc} approximations, either full or sparse, is to find a complete set of polynomials that are orthogonal to each other with respect to the input parameter \glspl{pdf}. 
In the case of arbitrary \glspl{pdf}, such polynomials can be numerically constructed \cite{soize2004, wan2006a, wan2006}.
As a result, despite the fact that according to a number of studies Lagrange interpolation methods enjoy an advantage over \gls{gpc}-based {methods} with respect to the error-cost ratio \cite{eldred2009a, eldred2009, ng2012}, they are abandoned whenever the input parameter \glspl{pdf} are not uniform, or normal, or closely related distributions, e.g., log-uniform or log-normal.
Consequently, the literature lacks comparisons between Lagrange interpolations and \gls{gpc}-based approximations for {other} probability measures.

This work investigates the use of Lagrange interpolation methods based on the recently introduced weighted Leja node sequences \cite{narayan2014} for the approximation of systems with arbitrarily distributed random parameters.
Weighted Leja sequences provide interpolation and quadrature nodes tailored to an arbitrary \gls{pdf} in a seamless way. 
Moreover, being nested by definition, they allow for  {the efficient construction of sparse grids,  thus enabling high-dimensional \gls{uq} studies}.
 {We note that this paper focuses on continuous \glspl{pdf} of arbitrary shapes and does not consider data-driven approaches which define a \gls{pdf} through its moments without any further assumptions about it \cite{ahlfeld2016, oladyshkin2012}}.

Due to a number of desirable properties, e.g., granularity, interpolation, and quadrature stability, and nestedness, weighted Leja nodes have been getting increasing attention in the context of approximation and \gls{uq} \cite{chkifa2014, ernst2019, georg2018, narayan2014, nobile2015, schillings2013}.
Nonetheless, with the exceptions of \cite{georg2018, loukrezis2019, loukrezis2019phd} where beta distributions are considered for the parameters, the use of Leja nodes has been limited to uniform and log-normal parameter distributions.
The contribution of this paper is exactly to start filling this gap with the application of weighted Leja interpolation for possibly high-dimensional systems with arbitrary parameter \glspl{pdf}.  

To that end, we present here the use of weighted Leja rules in the context of sparse Lagrange interpolation for approximation and \gls{uq} purposes in a comprehensive and self-consistent way. 
Section \ref{sec:materials} summarizes the main theory related to weighted Leja interpolation.
In Section~\ref{sec:stoch_param} we first present a general overview of the \gls{uq} setting in which Leja interpolation will be applied.
Then, in Section~\ref{sec:univar} we recall the basics of univariate Lagrange interpolation, both standard and hierarchical.
Unweighted and weighted Leja rules along with their properties are presented in Section~\ref{sec:leja}, while in in Section~\ref{sec:sparse_interp} the extension of Lagrange interpolation to higher dimensions, along with a dimension-adaptive, Leja-based scheme for the construction of sparse approximations and its use for \gls{uq} purposes, is presented.
Numerical experiments which verify the suitability of the Leja interpolation scheme for physical systems with arbitrarily distributed stochastic input parameters are presented in Section~\ref{sec:num_exp}.
In the same section, the adaptive weighted Leja interpolation \gls{uq} method is compared against a well-known adaptive \gls{gpc} algorithm based on \gls{lar}.
We conclude with a discussion on the overall findings of this work and on possibilities for further research in Section~\ref{sec:concl}.

\section{Theory and Methodology}
\label{sec:materials}
\subsection{Stochastic Parametric Models}
\label{sec:stoch_param}
Without loss of generality, let us assume that the behavior of the physical system under consideration can be modeled by a set of parametric \glspl{pde}.
The general form of the corresponding mathematical model is given by:
\begin{equation}
\label{eq:general}
\mathcal{D}_{\mathbf{y}}(u) = 0,
\end{equation}
where $\mathcal{D} = \mathcal{D}_{\mathbf{y}}(u)$ is a differential operator,  $\mathbf{y} \in \Xi \subset \mathbb{R}^N$ a $N$-dimensional parameter vector, $u = u_{\mathbf{y}}$ the parameter-dependent solution of Equation \eqref{eq:general},  {and $\Xi$ a measurable set called the image space}. %We added ``Equation'', please confirm.
Typically, we are interested in a {certain quantity, that is} the \gls{qoi}, which is computed by post-processing the solution $u = u_{\mathbf{y}}$.
For simplicity, we assume that the \gls{qoi} is a real scalar, however, the methodology discussed in Section~\ref{sec:sparse_interp} can be applied to complex or vector-valued \glspl{qoi} with only minor modifications.
We further assume that the functional dependence between the parameter vector and the \gls{qoi} is given by the map $g: \mathbf{y} \mapsto g\left(\mathbf{y}\right)$, which is assumed to be deterministic, i.e., the exact same output $g\left(\mathbf{y}\right)$ is observed each and every time the \gls{qoi} is evaluated for the same parameter vector $\mathbf{y}$.

Parametric problems in the form of Equation \eqref{eq:general} arise in both deterministic and stochastic settings. 
An example of the former case is, e.g., an optimization study, while an example of the latter is, e.g., a \gls{uq} study as considered in this work.
In the stochastic setting, the parameter vector $\mathbf{y}$ corresponds to a realization of a $N$-dimensional \gls{rv} $\mathbf{Y} = \left(Y_1, Y_2, \dots, Y_N\right)$, alternatively called a random vector.
Assuming that the \glspl{rv} $Y_n$, $n=1,\dots,N$, are mutually independent, the joint \gls{pdf} is given as a product of univariate \glspl{pdf}, such that $\varrho(\mathbf{y}) = \prod_{n=1}^N \varrho_n(y_n).$
While the case of dependent \glspl{rv} is not  {rigorously addressed here}, we mention the option of using generalized Nataf or Rosenblatt transformations to transform the dependent \glspl{rv} into independent \glspl{rv} \cite{feinberg2018, lebrun2009}.
The sparse interpolation method presented in Section~\ref{sec:sparse_interp} can then be applied to the transformed \glspl{rv}.

Due to its dependence on stochastic input parameters, the deterministic map $g$ results now in a random output. 
In other words, the input uncertainty propagates through the deterministic model and renders the \gls{qoi} uncertain.
The main goal of \gls{uq} studies is to quantify this output uncertainty by computing statistics, such as moments, event probabilities, or sensitivity metrics.
Denoting the expectation operator {by} $\mathbb{E}\left[\cdot\right]$, any statistic can be written in the general form:
\begin{equation}
\label{eq:stat_measure}
\mathbb{E}\left[\phi\left(g\right)\right] = \int_{\Xi} \phi\left(g\left(\mathbf{y}\right)\right)  \varrho\left(\mathbf{y}\right) \mathrm{d}\mathbf{y},
\end{equation}
where the functional $\phi$ corresponds to the {statistic of interest, e.g. $\phi\left(g\right) = g$ for the expected value of the \gls{qoi} and $\phi\left(g\right) = \left(g - \mathbb{E}\left[g\right]\right)^2$ for its variance.}

The integral of Equation  \eqref{eq:stat_measure} cannot in general be computed analytically. Therefore, one relies on numerical integration methods, such as \gls{mc} integration or quadrature.
Considering the case of \gls{mc} integration, if the original input-output map $g$ corresponds to an expensive numerical model,  {then it is desirable to replace it by} an inexpensive albeit accurate surrogate model $\widetilde{g} \approx g$.
In this work, this surrogate is obtained by interpolating $g$.
To capture accurately the behavior of $g(\mathbf{y})$, given that $\mathbf{y}$ are realizations of a random vector, the interpolation must be based on nodes tailored to the \gls{pdf} $\varrho(\mathbf{y})$.
The surrogate model can also be used for tasks other than computing statistics, e.g., for optimization.
Moreover, by integrating the interpolant, one can derive a quadrature scheme for the computation of Equation \eqref{eq:stat_measure}. Similar to the interpolation,  {nodes tailored to the \gls{pdf} are} crucial for the convergence of the quadrature.
Finally, in the case of a high-dimensional random vector $\mathbf{Y}$, sparse quadrature and interpolation schemes must be used in order to mitigate the effect of the curse of dimensionality. 
Such sparse schemes most commonly rely on nested univariate node sequences (see Section \ref{sec:univar}).

As pointed out in Section \ref{sec:intro}, there exist interpolation nodes that meet the aforementioned requirements for the well studied cases of uniform and normal \glspl{pdf}, e.g., Clenshaw--Curtis and Genz--Keister nodes, respectively \cite{clenshaw1960,genz1996}.
We note that only a finite number of Genz--Keister nodes {are} available, thus imposing an additional limitation on its achievable accuracy.
{However, constructing suitable interpolation and quadrature rules for more exotic or even arbitrary \glspl{pdf}, while in principle possible \cite{burkardt2014, sommariva2013}, remains a computational challenge.}
This problem is here circumvented with the use of weighted Leja nodes, discussed in Section \ref{sec:leja}, which are easy to compute, by definition nested, and tailored to a continuous \gls{pdf}.

\subsection{Univariate Interpolation Schemes}
\label{sec:univar}
Let us assume a univariate input-output map $g\left(y\right)$, where $y$ denotes a realization of a single \gls{rv} $Y$.
We let the non-negative integer $i \in \mathbb{Z}^{\geq 0}$ denote the interpolation level, $Z_i = \left\{y_{i,j}\right\}_{j=1}^{m(i)}$ the level-$i$ grid of interpolation nodes, and $m:\mathbb{Z}^{\geq 0} \rightarrow \mathbb{Z}^{+}$, $m\left(0\right) = 1$, a strictly monotonically increasing ``level-to-nodes'' function, relating the interpolation level to the grid size.
Based on these three elements, in the following subsections we define Lagrange and hierarchical interpolation rules, as well as related quadrature schemes. 

\subsubsection{Univariate Lagrange Interpolation}
\label{subsec:univar_lagrange}
The level-$i$ univariate Lagrange interpolation is defined as:
\begin{equation}
\label{eq:interp1d}
\mathcal{I}_{i} \left[g\right]\left(y\right) = \sum_{j=1}^{m(i)} g\left(y_{i,j}\right) l_{i, j}\left(y\right),
\end{equation}
where $l_{i, j}$ are univariate Lagrange polynomials, defined on the interpolation grid $Z_i$ as:
\begin{equation*}
\label{eq:lagrange1d}
l_{i, j}\left(y\right) = 
\prod_{k=1, k \neq j}^{m(i)} \frac{y-y_{i,k}}{y_{i, j} - y_{i,k}}.
\end{equation*}
The accuracy of the interpolation depends crucially on the choice of nodes, e.g., Gauss--Legendre and Gauss--Hermite quadrature nodes are popular choices for \glspl{rv} following uniform or normal distributions, respectively.
In particular, the error between $g$ and $\mathcal{I}_i\left[g\right]$ is related to the best polynomial approximation error $\epsilon_{\text{best}}$ through \cite{barthelmann2000, chkifa2014}:
\begin{equation*}
\label{eq:approx_error}
\left\lVert g - \mathcal{I}_i\left[g\right] \right\rVert_{L^{\infty}\left(\mathbb{P}_i\right)} \leq \left(1 + \mathbb{L}_i \right) \epsilon_{\text{best}},
\end{equation*}
where $\mathbb{P}_i$ denotes the space of Lagrange polynomials with degree up to $m(i) - 1$ and $\mathbb{L}_i$ the Lebesgue constant:
\begin{equation*}
\label{eq:Lebesgue}
\mathbb{L}_i = \mathrm{sup} \frac{\left\lVert \mathcal{I}_i\left[g\right] \right\rVert_{L^{\infty}\left(\mathbb{P}_i\right)}}{\left\lVert g \right\rVert_{L^{\infty}\left(\mathbb{P}_i\right)}}.
\end{equation*} 
Therefore, assuming that the interpolation error decreases for increasing interpolation levels, the interpolation nodes should be associated to a Lebesgue constant which grows at a slower rate.

\subsubsection{Hierarchical Univariate Interpolation}
\label{subsec:univar_hierarch}
We further demand that the interpolation nodes form sequences of nested grids for increasing interpolation levels, such that $Z_{i-1} \subset Z_{i}$, $\forall i >0$.
In that case, we may define a hierarchical interpolation scheme as follows.
We first define the operator $\Delta_i$ as the difference between two consecutive univariate interpolants, i.e.,
\begin{equation}
\label{eq:diff1d}
\Delta_i = \mathcal{I}_i - \mathcal{I}_{i-1},
\end{equation}
where $\mathcal{I}_{-1}\left[g\right](y) = 0$.
The interpolation in Equation  \eqref{eq:interp1d} can now be given as a sum of interpolant differences, such that:
\begin{equation}
\label{eq:hinterp1d_1}
\mathcal{I}_{i} \left[g\right]\left(y\right) = \sum_{k=0}^i \Delta_k\left[g\right]\left(y\right).
\end{equation}  
The hierarchy of formula in Equation  \eqref{eq:hinterp1d_1} is obvious if we consider two consecutive interpolation levels, $i-1$ and $i$, in which case the level-$i$ interpolation reads:
\begin{subequations}
	\label{eq:hinterp1d_2}
	\begin{align}
	\mathcal{I}_{i} \left[g\right]\left(y\right) &= \mathcal{I}_{i-1} \left[g\right]\left(y\right) + \sum_{j:y_{i,j} \in Z_i \setminus Z_{i-1}} \left( g\left(y_{i,j}\right) - \mathcal{I}_{i-1} \left[g\right]\left(y_{i,j}\right)\right) l_{i,j}\left(y\right) \\
	&= \mathcal{I}_{i-1} \left[g\right]\left(y\right) + \sum_{j:y_{i,j} \in Z_i \setminus Z_{i-1}} s_{i,j} l_{i,j}\left(y\right),
	\end{align}
\end{subequations}
where the interpolation coefficients $s_{i,j}$, given by: 
\begin{equation}
\label{eq:hs1d}
s_{i,j} = g\left(y_{i,j}\right) - \mathcal{I}_{i-1} \left[g\right]\left(y_{i,j}\right), \:\:\: y_{i,j} \in Z_i \setminus Z_{i-1},
\end{equation}
are called hierarchical surpluses. 
In the particular case of $m(i) = i + 1$, which can in fact be accomplished using the Leja sequences discussed in Section~\ref{sec:leja}, a single grid point is added for each new interpolation level. Then, the index $j$ can be dropped and Equation \eqref{eq:hinterp1d_2} is simplified to  {$\mathcal{I}_{i} \left[g\right]\left(y\right) = \mathcal{I}_{i-1} \left[g\right]\left(y\right) + s_{i} l_{i}\left(y\right)$.}
%\begin{align}
%\label{eq:hinterp1d_3}
%\mathcal{I}_{i} \left[g\right]\left(y\right) &= \mathcal{I}_{i-1} \left[g\right]\left(y\right) +  \left( g\left(y_{i}\right) - \mathcal{I}_{i-1} \left[g\right]\left(y_{i}\right)\right) l_{i}\left(y\right) = \mathcal{I}_{i-1} \left[g\right]\left(y\right) + s_{i} l_{i}\left(y\right).
%\end{align}

The obvious advantage of using nested grids and the hierarchical Equations  \eqref{eq:hinterp1d_1} or \eqref{eq:hinterp1d_2}, is that at each new level $i$ the \gls{qoi} must be evaluated only at the new nodes $y_{i,j} \in Z_i \setminus Z_{i-1}$.
Moreover, univariate hierarchical interpolation rules constitute the backbone of the sparse adaptive multivariate interpolation algorithms discussed in Section~\ref{sec:sparse_interp}.

\subsubsection{Interpolatory Univariate Quadrature}
\label{subsec:univar_quad}
Considering the \gls{uq} goal of estimating specific statistics given by  Equation \eqref{eq:stat_measure}, we may use an available interpolation $\mathcal{I}_i\left[g\right]$ to derive an appropriate quadrature scheme.
For example, assuming an univariate interpolation given by:
\begin{equation*}
\mathcal{I}_{i} \left[g\right]\left(y\right) = \sum_{j=1}^{m(i)} s_{i,j} l_{i, j}\left(y\right),
\end{equation*}
the expected value of the \gls{qoi} can be estimated as:
\begin{equation}
\label{eq:mean_interp_1d} 
\mathbb{E}\left[g\right] \approx \mathbb{E}\left[\mathcal{I}_i\left[g\right]\right] = \sum_{j=1}^{m(i)} s_{i,j} w_{i,j},
\end{equation}
where $w_{i,j}$ are quadrature weights given by:
\begin{equation*}
w_{i,j} = \mathbb{E}\left[l_{i,j}\right] =  \int_{\Xi} l_{i,j}\left(y\right) \varrho(y) \mathrm{d}y.
\end{equation*}
Quadrature schemes for the estimation of statistics other than the expected value can be derived in a similar way, by applying the expectation operator along with the associated functional $\phi$.

\subsection{Leja Interpolation Nodes}
\label{sec:leja}
As pointed out in Sections~\ref{sec:stoch_param} and \ref{sec:univar}, the univariate interpolation nodes employed in Lagrange interpolation schemes must satisfy certain requirements.
First, for the interpolation to become more accurate with increasing levels, the Lebesgue constant associated with the interpolation nodes must increase at a rate lower than the decrease in the polynomial approximation error \cite{barthelmann2000, chkifa2014}.
Second, the choice of nodes must result in accurate quadrature schemes, similar to Equation \eqref{eq:mean_interp_1d}, to be used in the computation of statistics.
Third, interpolation nodes which form nested grids can be used to derive hierarchical interpolation schemes, which result in reduced computations per added interpolation level.
Additionally, as discussed in Section~\ref{sec:sparse_interp}, sparse interpolations depend crucially on hierarchical univariate interpolation rules.
Finally and most importantly in the context of this work, the interpolation must be constructed with respect to arbitrary weight functions, i.e., input \glspl{pdf}.
All requirements are addressed here by using weighted Leja sequences \cite{narayan2014}, discussed in the following.

\subsubsection{Unweighted Leja Nodes}
\label{subsec:leja_unweighted}
A sequence of standard, unweighted Leja nodes \cite{leja1957} $\left\{y_j\right\}_{j \geq 0}$, $y_j \in \left[-1,1\right]$, is defined by solving the optimization problem:
\begin{equation}
\label{eq:lejaopt_unweighted}
y_{j} = \argmax_{y \in \left[-1, 1\right]} \prod_{k=0}^{j-1}\left|y - y_k\right|,
\end{equation}
for each node $y_j$, $j \geq 1$.
The initial node $y_0 \in \left[-1,1\right]$ can be chosen arbitrarily, therefore Leja sequences are not unique.
Using simple scaling rules, the unweighted Leja sequences defined in Equation \eqref{eq:lejaopt_unweighted} are directly applicable to Lagrange interpolation involving constant weight functions, equivalently, uniform \glspl{pdf} in the \gls{uq} context.
The extension to non-uniform \glspl{pdf} is discussed in Section~\ref{subsec:leja_weighted}.
With respect to the remaining node requirements, it has been shown that the Lebesgue constant associated with unweighted Leja nodes increases subexponentially \cite{taylor2010}. 
Regarding Leja-based quadrature rules, an extensive discussion  can be found in \cite{narayan2014}. 
Considering uniform \glspl{pdf}, a number of works \cite{loukrezis2019, narayan2014, schillings2013} show that Leja-based quadrature schemes are sufficiently accurate, albeit suboptimal compared to more standard choices, e.g., based on Clenshaw--Curtis nodes.
Finally, the optimization problem of Equation \eqref{eq:lejaopt_unweighted} results in nested node sequences irrespective of the employed level-to-nodes function $m(i)$. 
This is an additional advantage of using Leja nodes, since the nested node sequences can be as granular as the user wishes. 
In this work, we use $m(i) = i+1$ to get the minimum of one extra node per interpolation level, i.e., $\# \left(Z_i \setminus Z_{i-1}\right) = 1$.
Another popular choice is to use $m(i) = 2i + 1$ \cite{schillings2013}, which is typically employed for the construction of symmetric Leja sequences, where $y_{0} = 0$, $y_{1} = 1$, $y_j$ is given as in Equation  \eqref{eq:lejaopt_unweighted} for odd $j>1$, and $y_j = -y_{j-1}$ for even $j$.
%\begin{subequations}
%	\label{eq:symmetric-Leja}
%	\begin{align}
%	y_{j} &= \argmax_{y \in \left[-1, 1\right]} \prod_{k=0}^{j-1}\left|y - y_k\right|, &&j=3,5,7,\dots, \\
%	y_{j} & = -y_{j-1}, &&j=2,4,6,\dots,
%	\end{align}
%\end{subequations}
%as well as for $\mathcal{R}$-Leja sequences \cite{calvi2011}, defined as projections on $\left[-1,1\right]$ of Leja sequences computed on the complex unit disc, such that
%\begin{subequations}
%	\label{eq:R-Leja}
%	\begin{align}
%	y_{j} &= \mathcal{R}\left(\argmax_{\left|z\right| \leq 1} \prod_{k=0}^{j-1}\left|z - z_k\right|\right), &&j=3,5,7,\dots, \\
%	y_{j} & = -y_{j-1}, &&j=2,4,6,\dots,
%	\end{align}
%\end{subequations}
%where $y_{0} = 0$, $y_{1} = 1$, and $\mathcal{R}\left(z\right)$ denotes the real part of a complex number $z$.
In all cases, Leja sequences are more granular compared to other nested node sequences, e.g., the level-to-nodes function $m(i) = 2^i + 1$ must be used for nested Clenshaw--Curtis nodes.

\subsubsection{Weighted Leja Nodes}
\label{subsec:leja_weighted}
In more general \gls{uq} settings, e.g., for the arbitrary input distributions considered in this work, the definition of Leja sequences must be adjusted according to the input \gls{pdf} $\varrho\left(y\right)$, which acts as a weight function.
To that end, we employ the definition of weighted Leja sequences given in \cite{narayan2014}.
Using this formulation, weighted Leja sequences are constructed by incorporating the \gls{pdf} $\varrho\left(y\right)$ into the optimization problem of Equation  \eqref{eq:lejaopt_unweighted}, which is transformed into
\begin{equation}
\label{eq:lejaopt_weighted}
y_{j} = \argmax_{y \in \Xi} \sqrt{\varrho\left(y\right)} \prod_{k=0}^{j-1}\left|y - y_k\right|.
\end{equation}
We note that other formulations exist as well \cite{demarchi2004, jantsch2019}, however, the form given in Equation  \eqref{eq:lejaopt_weighted} is preferred due to the fact that the resulting Leja nodes are asymptotically distributed as weighted Gauss quadrature nodes \cite{narayan2014}. 
While this is a promising property, it is not necessarily sufficient to produce an accurate interpolation rule.
In some cases, e.g., for weight functions resembling the Gaussian \gls{pdf} \cite{jantsch2019}, it can be shown that weighted Leja nodes are also associated with a subexponentially growing Lebesgue constant. 
While more general theoretical results are currently not  available, weighted Leja interpolation has been found to perform very  {well} in practice \cite{ernst2019, loukrezis2019}, as also demonstrated {in this paper} by the results of the numerical experiments presented in Section~\ref{sec:num_exp}.
The other beneficial properties of Leja sequences, i.e., granularity and nestedness, are preserved in the weighted case. 
Moreover, since the interpolation rule is now tailored to the specific weight function that is the \gls{pdf}, the corresponding quadrature rules will be tailored to that \gls{pdf} as well.
For a more detailed discussion on the properties of weighted Leja sequences, see \cite{narayan2014}.

\subsection{Sparse Adaptive Leja Interpolation}
\label{sec:sparse_interp}
The multivariate interpolation is based on suitable combinations of the univariate interpolation rules discussed in Section~\ref{sec:univar}.
We denote with $\mathbf{i} = \left(i_1, i_2, \dots, i_N\right)$ a multi-index of interpolation levels per input \gls{rv}, with $\Lambda$ a multi-index set, and with $\Delta_{\mathbf{i}}$ the tensor-product difference operator given by:
\begin{equation*}
\label{eq:diffNd}
\Delta_{\mathbf{i}} = \Delta_{1,i_1} \otimes \Delta_{2,i_2} \otimes \cdots \otimes \Delta_{N,i_N},
\end{equation*}
where the univariate difference operators $\Delta_{n, i_n}$ are defined as in Equation \eqref{eq:diff1d}, to be recalled in the following subsections.
In Section~\ref{subsec:gensparse}, we discuss interpolation on generalized sparse grids.
A dimension-adaptive interpolation scheme based on Leja sequences is presented in Section~\ref{subsec:dimadapt}.
Post-processing multivariate interpolations for \gls{uq} purposes is discussed in Section~\ref{subsec:postproc}.

\subsubsection{Generalized Sparse Grid Interpolation}
\label{subsec:gensparse}
Given a multi-index set $\Lambda$, the interpolation-based approximation of the multivariate input-output map $g(\mathbf{y})$ reads:
\begin{equation}
\label{eq:sparseint}
\mathcal{I}_{\Lambda}\left[g\right]\left(\mathbf{y}\right) = \sum_{\mathbf{i} \in \Lambda} \Delta_{\mathbf{i}}\left[g\right]\left(\mathbf{y}\right).
\end{equation}
We note that Equation  \eqref{eq:sparseint} is not in general interpolatory, i.e., the interpolation property $\mathcal{I}_{\Lambda}\left[g\right]\left(\mathbf{y}_*\right) = g\left(\mathbf{y}_*\right)$, where $\mathbf{y}_*$ is a multivariate interpolation node, does not necessarily hold. 
The interpolation property holds only if $\Lambda$ is a downward-closed set and the underlying univariate interpolation rules are based on nested nodes \cite{barthelmann2000}.
Downward-closed sets, also known as monotone or lower sets, satisfy the property:
\begin{equation}
\label{eq:monotonicity}
\forall \mathbf{i} \in \Lambda \Rightarrow \mathbf{i}-\mathbf{e}_n \in \Lambda, \forall n = 1,2,\dots,N, \: \text{with} \:\: i_n > 0,
\end{equation}
where $\mathbf{e}_n$ denotes the unit vector in the $n$th dimension.
If, additionally, the multi-index set $\Lambda$ satisfies some sparsity constraint, then Equation \eqref{eq:sparseint} is a sparse interpolation formula, while the corresponding sparse grid is given by:
\begin{equation*}
\label{eq:gridNd}
Z_{\Lambda} = \bigcup_{\mathbf{i} \in \Lambda} Z_{\mathbf{i}} = \bigcup_{\mathbf{i} \in \Lambda} \left(Z_{1,i_1} \times Z_{2,i_2} \times \cdots \times Z_{N,i_N}\right).
\end{equation*}
A commonly used sparsity constraint is: 
\begin{equation*}
\label{eq:isosparse} 
\Lambda = \left\{\mathbf{i} : \sum_{n=1}^N i_n \leq i_{\max}\right\},
\end{equation*}
resulting in the so-called isotropic sparse grids \cite{babuska2010, barthelmann2000}. 
Compared to full isotropic tensor grids, i.e., those given by:
\begin{equation*}
\label{eq:tensorproduct} 
\Lambda = \left\{\mathbf{i} : \max_{n=1,\dots,N} i_n \leq i_{\max}\right\},
\end{equation*}
the complexity of which is $\mathcal{O}\left(m \left(i_{\max}\right)^N\right)$, isotropic sparse grids delay significantly the curse of dimensionality, since their complexity is reduced to $\mathcal{O}\left(m \left(i_{\max}\right) \left(\log m \left(i_{\max}\right)\right)^{N-1}\right)$ \cite{bungartz2004}.
Nonetheless, isotropic sparse grids still grow exponentially with respect to the number of parameters.

In most practical cases, the influence of the input parameters upon the \gls{qoi} is anisotropic. 
For example, the \gls{qoi} may be more sensitive to some parameters compared to others, in the sense that variations in the first set of parameters result in comparably greater variations in the \gls{qoi}. 
Accordingly, the functional dependence of the \gls{qoi} on  {a certain} parameter might be significantly more difficult to approximate compared to another, e.g., a highly nonlinear versus a linear relation.

The interpolation of Equation  \eqref{eq:sparseint} can be constructed such that this anisotropy is represented in its terms, equivalently in the multi-indices comprising the set $\Lambda$ \cite{chkifa2014, nobile2008a}.
This is accomplished via an anisotropic refinement of the parameter space according to the contribution of each parameter to the interpolation accuracy. 
Compared to isotropic schemes, anisotropic interpolations typically result in great computational savings for the same approximation accuracy, or, equivalently, in greater interpolation accuracy for an equal cost.
Since, in most cases, the parameter anisotropy cannot be a priori estimated in an accurate way \cite{beck2012}, anisotropic interpolations are usually constructed using greedy, adaptive algorithms \cite{gerstner2003, klimke2005, stoyanov2016}. 
For the particular case of Leja nodes, such an algorithm is discussed in Section~\ref{subsec:dimadapt}.
Leja-based, adaptive, anisotropic sparse-grid interpolation algorithms are also available in \cite{chkifa2014, narayan2014}.

\subsubsection{Adaptive Anisotropic Leja Interpolation}
\label{subsec:dimadapt}
We will base the adaptive construction of anisotropic sparse interpolations on the dimension-adaptive algorithm first presented in \cite{gerstner2003} for quadrature purposes.
Variants of this algorithm appear in a number of later works \cite{chkifa2014, klimke2005, loukrezis2019, narayan2014, schillings2013}.
In this work, we assume the all underlying univariate interpolation rules employ Leja nodes, presented in Section~\ref{sec:leja}, along with the level-to-nodes function $m(i) = i+1$.
Our approach is depicted in Algorithm~\ref{algo:dimadapt}, while a detailed presentation follows.
\RestyleAlgo{boxruled}

\begin{algorithm}[b!]
\caption{Dimension-adaptive Leja interpolation.} %please change algorithm 1 to tabular form.
% authors' comment: we do not undertand the reviewers comment. the algorithm is already in the format given in the example
	\SetAlgoLined
	\KwData{input-output map $g\left(\mathbf{y}\right)$, initial downward-closed multi-index set $\Lambda^{\text{init}}$, tolerance $\epsilon$, simulation budget $B$.}
	\KwResult{multi-index set $\Lambda$, sparse interpolation $\mathcal{I}_{\Lambda}\left[g\right]$, sparse grid $Z_{\Lambda}$.}
	$ \Lambda \leftarrow \Lambda^{\text{init}}$.\\
	\While{\textsc{True}}{
		Compute admissible set: $\Lambda^{\text{adm}} = \left\{\mathbf{i} : \mathbf{i} \not \in \Lambda \:\: \text{and} \:\: \Lambda \cup \left\{\mathbf{i}\right\} \:\: \text{is downward-closed}\right\}$. \\
		Compute admissible Leja nodes: $\mathbf{y}_{\mathbf{i}} = Z_{\mathbf{i}} \setminus Z_{\Lambda}$, $\forall \mathbf{i} \in \Lambda^{\text{adm}}$.\\
		Compute hierarchical surpluses: $ s_{\mathbf{i}} = g\left(\mathbf{y}_{\mathbf{i}}\right) - \mathcal{I}_{\Lambda}\left[g\right]\left(\mathbf{y}_{\mathbf{i}}\right)$, $\forall \mathbf{i} \in \Lambda^{\text{adm}}$. \\
		Compute contribution indicators: $ \eta_{\mathbf{i}} = \left|s_{\mathbf{i}}\right|$, $\forall \mathbf{i} \in \Lambda^{\text{adm}}$. \\
		Compute costs: $ C = \#Z_{\Lambda} + \#Z_{\Lambda^{\text{adm}}}$. \\
		Compute total contribution of the admissible set:  $\eta^{\text{tot}} = \sum_{\mathbf{i} \in \Lambda^{\text{adm}}} \eta_{\mathbf{i}}$. \\
		\If{$C \geq B$ \textsc{or} $\eta^{\mathrm{tot}} \leq \epsilon$}{Exit while-loop.}
		Find new multi-index: $ \mathbf{i}^*  = \argmax_{\mathbf{i} \in \Lambda^{\text{adm}}}\eta_{\mathbf{i}}$. \\
		Update multi-index set: $ \Lambda \leftarrow \Lambda \cup \left\{\mathbf{i}^*\right\}$. \\
	}
	Construct final multi-index set $\Lambda \leftarrow \Lambda \cup \Lambda^{\text{adm}}$, sparse grid $Z_\Lambda$, and approximation $\mathcal{I}_{\Lambda}\left[g\right]$.
	\label{algo:dimadapt}
\end{algorithm}

Let us assume that a sparse interpolation $\mathcal{I}_{\Lambda}$, based on a downward-closed multi-index set $\Lambda$ and on univariate Leja sequences, is given as in Equation  \eqref{eq:sparseint}.
We define the set of admissible multi-indices, i.e., those that satisfy the downward-closedness property of Equation  \eqref{eq:monotonicity} if added to $\Lambda$, as: 
\begin{equation*}
\label{eq:admset}
\Lambda^{\text{adm}} = \left\{\mathbf{i} : \mathbf{i} \not \in \Lambda \:\: \text{and} \:\: \Lambda \cup \left\{\mathbf{i}\right\} \:\: \text{is downward-closed}\right\}.
\end{equation*}
Due to the fact that all univariate Leja rules employ the level-to-nodes function $m(i) = i+1$, each multi-index $\mathbf{i} \in \Lambda^{\text{adm}}$ is associated with a single new interpolation node $\mathbf{y}_{\mathbf{i}} = Z_{\mathbf{i}} \setminus Z_{\Lambda}$.
Therefore, should an admissible multi-index $\mathbf{i} \in \Lambda^{\text{adm}}$ be added to $\Lambda$, we obtain the hierarchical interpolation scheme:
\begin{equation*}
\label{eq:hierarch_Nd}
\mathcal{I}_{\Lambda \cup \left\{ \mathbf{i} \right\} }\left[g\right]\left(\mathbf{y}\right) = \mathcal{I}_{\Lambda}\left[g\right]\left(\mathbf{y}\right) +  s_{\mathbf{i}} L_{\mathbf{i}}\left(\mathbf{y}\right),
\end{equation*}
where, similarly to Equation  \eqref{eq:hs1d}, the hierarchical surpluses $s_{\mathbf{i}}$ are given by:
\begin{equation*}
\label{eq:hsNd}
s_{\mathbf{i}} = g\left(\mathbf{y}_{\mathbf{i}}\right) - \mathcal{I}_{\Lambda} \left[g\right] \left(\mathbf{y}_{\mathbf{i}}\right), \quad \mathbf{i} \in \Lambda^{\text{adm}},
\end{equation*}
and $L_{\mathbf{i}}$ are multivariate Lagrange polynomials given by:
\begin{equation*}
L_{\mathbf{i}}(\mathbf{y}) = \prod_{n=1}^N l_{n, i_n}(y_n), 
\end{equation*}
where, due to the use of nested Leja sequences with the level-to-nodes function $m(i)=i+1$, the univariate Lagrange polynomials can now be defined hierarchically \cite{chkifa2014}, such that,
\begin{equation*}
\label{eq:lagrange1d_hierarch}
l_{n,i_n}\left(y_n\right) = 
\prod_{k=0}^{i_n-1} \frac{y_n-y_{n,k}}{y_{n,i_n} - y_{n,k}}.
\end{equation*}

Similar to \cite{gerstner2003}, we expect that coefficients of small magnitudes correspond to terms with relatively small contributions to the available interpolation, which can thus be omitted.
Accordingly, terms whose coefficients have large magnitudes are expected to enhance the accuracy of the interpolation and should therefore be added.
Based on this argument, we define for each multi-index $\mathbf{i} \in \Lambda^{\text{adm}}$ the contribution indicator:
\begin{equation*}
\label{eq:local_err_ind}
\eta_{\mathbf{i}} = \left|s_{\mathbf{i}}\right|.
\end{equation*}
Naturally, the multi-index set $\Lambda$ should be enriched with the admissible multi-index:
\begin{equation*}
\mathbf{i}^* = \argmax_{\mathbf{i} \in \Lambda^{\text{adm}}} \eta_{\mathbf{i}},
\end{equation*}
that is the index corresponding to the maximum contribution.
This procedure can be continued iteratively until a simulation budget $B$ is reached, or until the total contribution of the admissible set is below a tolerance $\epsilon$.
At any given step of the algorithm, the total number of simulations, equivalently, model evaluations or function calls, is equal to $\#Z_{\Lambda \cup \Lambda^{\text{adm}}}$. 
After the termination of the iterations, the final interpolation is constructed using all multi-indices in $\Lambda \cup \Lambda^{\text{adm}}$, since the hierarchical surpluses corresponding to the admissible nodes have already been computed.

\subsubsection{Post-Processing}
\label{subsec:postproc}
Once available, the sparse interpolation can be used as an inexpensive surrogate model that can replace the original model in computationally demanding tasks, such as sampling-based estimation of statistics.
Specific statistics can be computed directly after the interpolation terms.
One possibility is to first transform the Lagrange basis into a \gls{gpc} basis \cite{buzzard2013} and then use the well-known formulas of \gls{gpc} approximations for the estimation of moments or sensitivity indices \cite{blatman2010, blatman2011}.
Another approach is to apply the expectation operator along with the functional corresponding to the statistic directly upon the interpolation, as to derive an appropriate quadrature scheme, similar to the univariate case presented in Section~\ref{subsec:univar_quad}.
For example, assuming a Leja-based interpolation given by:
\begin{equation*}
\mathcal{I}_{\Lambda}\left[g\right]\left(\mathbf{y}\right) = \sum_{\mathbf{i} \in \Lambda} s_{\mathbf{i}} L_{\mathbf{i}}\left(\mathbf{y}\right),
\end{equation*}
the expected value of the \gls{qoi} can be estimated as:
\begin{equation}
\label{eq:mean_interp} 
\mathbb{E}\left[g\right] \approx \mathbb{E}\left[\mathcal{I}_{\Lambda}\left[g\right]\right] = \sum_{\mathbf{i} \in \Lambda} s_{\mathbf{i}} w_{\mathbf{i}},
\end{equation}
where the quadrature weights $w_{\mathbf{i}}$ are given as products of the respective univariate weights, i.e., 
\begin{equation*}
w_{\mathbf{i}} = \mathbb{E}\left[L_{\mathbf{i}}\right] =  \mathbb{E}\left[\prod_{n=1}^N l_{n, i_n}\right] = \prod_{n=1}^N \mathbb{E}\left[l_{n, i_n}\right] = \prod_{n=1}^N w_{n, i_n}.
\end{equation*}
Clearly, analytical formulas similar to Equation  \eqref{eq:mean_interp} can be derived for other statistics as well.

\section{Results}
\label{sec:num_exp}
\begin{figure}[H]
	\begin{tabular}[b]{c}
		\begin{tikzpicture} 
		\begin{axis}[width=0.5\textwidth, height=0.2\textwidth, xlabel=$y$, yticklabels={ }, xtick={0,0.5,1,1.5,2,2.5,3}, ] 
		\addplot[color=black,mark=*, only marks] coordinates {(4.297284620882929390e-06,0) (1.178995324029824127e-01,0) (3.168822439334964547e-01,0) (7.911577423902187434e-01,0) (1.176697435663403013e+00,0) (1.530916320376070949e+00,0) (1.934103072417560298e+00,0) (2.254221928522588403e+00,0) (2.684622798199928440e+00,0) (2.999996184893634066e+00,0)};
		\node at (axis cs:7.91157742e-01, 0.5) {0};
		\node at (axis cs:4.29728462e-06, 0.5) {1};
		\node at (axis cs:2.25422193e+00, 0.5) {2};
		\node at (axis cs:2.99999618e+00, 0.5) {3};
		\node at (axis cs:3.16882244e-01, 0.5) {4};
		\node at (axis cs:1.53091632e+00, 0.5) {5};
		\node at (axis cs:2.68462280e+00, 0.5) {6};
		\node at (axis cs:1.17899532e-01, 0.5) {7};
		\node at (axis cs:1.17669744e+00, 0.5) {8};
		\node at (axis cs:1.93410307e+00, 0.5) {9};
		\end{axis}
		\end{tikzpicture} \\
	({\bf{a}})
		%{\small (a) Univariate Leja nodes for $\mathcal{TN}\left(0,1,0,3\right)$.}
	\end{tabular}
	\hfill
	\begin{tabular}[b]{c}
		\begin{tikzpicture} 
		\begin{axis}[width=0.5\textwidth, height=0.2\textwidth, xlabel=$y$, yticklabels={ }, xtick={-2,-1,0,1,2,3},] 
		\addplot[color=black,mark=*, only marks] coordinates {(-1.999994107460605575e+00,0) (-1.769663700752260871e+00,0) (-1.389049730568840513e+00,0) (-7.124069122324625525e-01,0) (5.078323883960644824e-02,0) (8.735853143442018354e-01,0) (1.340332778815786607e+00,0) (1.784124487642308576e+00,0) (2.494524891576599845e+00,0) (2.999996130165510255e+00,0)};
		\node at (axis cs:0.05078324, 0.5) {0};
		\node at (axis cs:-1.38904973, 0.5) {1};
		\node at (axis cs:1.78412449, 0.5) {2};
		\node at (axis cs:-1.99999411, 0.5) {3};
		\node at (axis cs:2.99999613, 0.5) {4};
		\node at (axis cs:0.87358531, 0.5) {5};
		\node at (axis cs:-0.71240691, 0.5) {6};
		\node at (axis cs:2.49452489, 0.5) {7};
		\node at (axis cs:-1.7696637, 0.5) {8};
		\node at (axis cs:1.34033278, 0.5) {9};
		\end{axis}
		\end{tikzpicture} \\
		({\bf{b}})
		%{\small (b) Univariate Leja nodes for $\mathcal{TN}\left(0,1,-2,3\right)$.}
	\end{tabular}
	\\
	\vspace{1em}
	\begin{tabular}[b]{c}
		\begin{tikzpicture} 
		\begin{axis}[width=0.5\textwidth, height=0.2\textwidth, xlabel=$y$, yticklabels={ },] 
		\addplot[color=black,mark=*, only marks] coordinates {(-3.816900448682774449e+00,0) (-4.740760134152830485e-01,0) (5.308862651729214122e+00,0) (1.085360141806010503e+01,0) (1.919933261323517826e+01,0) (3.503371745275804017e+01,0) (4.555214368967180150e+01,0) (5.844480397479604505e+01,0) (8.385426160439362775e+01,0) (9.510339834941498793e+01,0)};
		\node at (axis cs:5.30886265, 0.5) {0};
		\node at (axis cs:-0.47407601, 0.5) {1};
		\node at (axis cs:19.19933261, 0.5) {2};
		\node at (axis cs:35.03371745, 0.5) {3};
		\node at (axis cs:10.85360142, 0.5) {4};
		\node at (axis cs:58.44480398, 0.5) {5};
		\node at (axis cs:-3.81690045, 0.5) {6};
		\node at (axis cs:83.8542616, 0.5) {7};
		\node at (axis cs:45.55214369, 0.5) {8};
		\node at (axis cs:95.10339835, 0.5) {9};
		\end{axis}
		\end{tikzpicture} \\
		({\bf{c}})
		%{\small (c) Univariate Leja nodes for $\mathcal{G}\left(3,4\right)$.}
	\end{tabular}
	\hfill 
	\begin{tabular}[b]{c}
		\begin{tikzpicture} 
		\begin{axis}[width=0.5\textwidth, height=0.2\textwidth, xlabel=$y$, yticklabels={ },] 
		\addplot[color=black,mark=*, only marks] coordinates {(-2.908450120297196850e+00,0) (-1.237039676890177464e+00,0) (1.654431325864607061e+00,0) (4.426800301870607690e+00,0) (8.599665787408531159e+00,0) (1.651685622669509712e+01,0) (2.177607028486232821e+01,0) (2.822240110288359460e+01,0) (4.092712979776482030e+01,0) (4.655169377958035426e+01,0)};
		\node at (axis cs:1.65443133, 0.5) {0};
		\node at (axis cs:-1.23703968, 0.5) {1};
		\node at (axis cs:8.59966579, 0.5) {2};
		\node at (axis cs:16.51685623 , 0.5) {3};
		\node at (axis cs:4.4268003, 0.5) {4};
		\node at (axis cs:28.2224011, 0.5) {5};
		\node at (axis cs:-2.90845012, 0.5) {6};
		\node at (axis cs:40.9271298, 0.5) {7};
		\node at (axis cs:21.77607028, 0.5) {8};
		\node at (axis cs:46.55169378, 0.5) {9};
		\end{axis}
		\end{tikzpicture} \\
		({\bf{d}})
		%{\small (d) Univariate Leja nodes for $\mathcal{G}\left(0.5,2\right)$.}
	\end{tabular}
	\caption{ {A total of 10 first univariate Leja nodes with respect to different truncated normal and Gumbel probability distributions. The truncated normal distributions are denoted with $\mathcal{TN}\left(\mu,\sigma^2, l, u\right)$, where $\mu$ and $\sigma^2$ refer to the mean value and the variance of a normal distribution, while $l$ and $u$ are the lower and upper truncation limits. The Gumbel distributions are denoted with $\mathcal{G}\left(\ell, \beta\right)$, where $\ell$ is the location parameter and $\beta$ the scaling parameter. The annotation above each node denotes the order of appearance in the Leja sequence, where ``0'' refers to the initial node: ({\bf{a}}) $\mathcal{TN}\left(\mu=0,\sigma^2=1,l=0,u=3\right)$, ({\bf{b}}) $\mathcal{TN}\left(\mu=0,\sigma^2=1,l=-2,u=3\right)$, ({\bf{c}}) $\mathcal{G}\left(\ell=3, \beta=4\right)$,  (\textbf{b}) $\mathcal{G}\left(\ell=0.5, \beta=2\right)$.}}
	\label{fig:leja_nodes}
\end{figure}
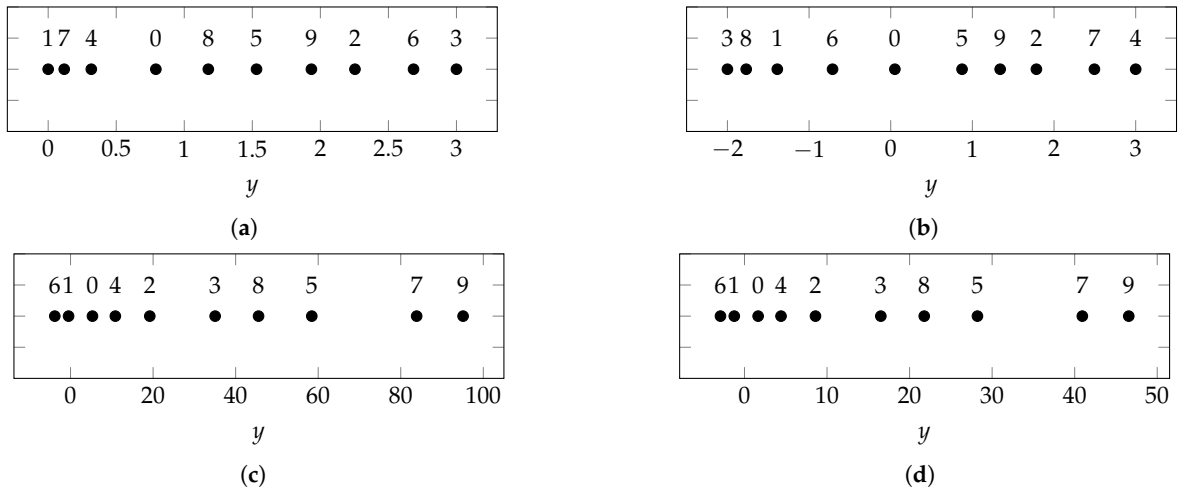

In this section, we consider models with relatively atypical input distributions, for which the application of interpolation-based \gls{uq} methods is not straightforward.
In particular, we consider input \glspl{rv} following truncated normal and Gumbel distributions, the latter also known as the generalized extreme value distribution of type 1 or the log-Weibull distribution. 
We denote the truncated normal distribution with $\mathcal{TN}\left(\mu,\sigma^2, l, u\right)$, where $\mu$ and $\sigma^2$ refer to the mean value and the variance of a normal distribution, while $l$ and $u$ are the lower and upper truncation limits. 
The Gumbel distribution is denoted with $\mathcal{G}\left(\ell, \beta\right)$, where $\ell$ is the location parameter and $\beta$ the scaling parameter.
Dedicated nested interpolation nodes are not readily available with respect to the corresponding input \glspl{pdf}, therefore, interpolation-based \gls{uq} methods have not been considered so far in the literature for such random inputs.
Here, we compute Leja nodes tailored to truncated normal and Gumbel distributions simply by applying Equation \eqref{eq:lejaopt_weighted} with weight functions that are equal to the corresponding \glspl{pdf}.
Exemplarily, Figure~\ref{fig:leja_nodes} shows the  first 10 computed Leja points with respect to two truncated normal and two Gumbel distributions, each with different distribution parameters.

Based on the numerical results presented in this section, we demonstrate that weighted Leja interpolation can reliably be used for such atypical input \glspl{pdf}, in terms of both approximation and \gls{uq}.
For further verification, we compare the performance of the weighted Leja interpolation against \gls{gpc} approximations with numerically constructed orthogonal polynomials \cite{soize2004}.
All models feature multiple input parameters, hence, we use adaptive algorithms resulting in sparse approximations.
The dimension-adaptive Algorithm~\ref{algo:dimadapt} based on weighted Leja nodes is implemented in our in-house developed DALI (Dimension-Adaptive Leja Interpolation) software \cite{loukrezis2019, loukrezis2019b}.
The sparse \gls{gpc} approximations are constructed with a well-established, degree-adaptive algorithm based on \glsfirst{lar} \cite{blatman2011} and implemented in the UQLab \cite{marelli2015} software.
Quasi-random experimental designs based on Sobol sequences are used in the \gls{lar}-\gls{gpc} approach, while numerically constructed orthogonal polynomials are used to tackle the input \glspl{pdf} \cite{soize2004}.

\subsection{Error Metrics}
\label{subsec:errors}
Let us assume that an interpolation or \gls{gpc}-based polynomial approximation $\widetilde{g} \approx g$ is available.
The following errors are used to estimate the performance of $\widetilde{g}$ for approximation and \gls{uq} purposes. 

First, we want to estimate the performance of $\widetilde{g}$ in terms of approximation accuracy, i.e., its suitability to act as a surrogate model and reliably replace the input-output map $g$, equivalently, the original computational model.
To that end, we use a validation sample $\left\{\mathbf{y}_q\right\}_{q=1}^Q$, which is randomly drawn from the joint input \gls{pdf} $\varrho(\mathbf{y})$ and compute the \gls{rms} validation error:
\begin{equation}
\label{eq:cverr_rms} 
\epsilon_{\text{RMS}} = \sqrt{\frac{1}{Q} \sum_{q=1}^Q \left(\widetilde{g}\left(\mathbf{y}_q\right) - g\left(\mathbf{y}_q\right)\right)^2}.
\end{equation}
From an interpretation point of view, $\epsilon_{\text{RMS}}$ indicates the average (expected) approximation accuracy of $\widetilde{g}$, while at the same time being sensitive to outliers by penalizing large errors.
In all numerical experiments, the size of the validation sample is set to $Q=10^5$.

Secondly, we want to estimate the performance of $\widetilde{g}$ in terms of \gls{uq} quality, i.e., estimate the accuracy of the statistics values provided by post-processing the approximation.
For that purpose, we use the expected value of the \gls{qoi} as a representative statistic and compute the relative error:
\begin{equation}
\label{eq:relerr_mean}
\epsilon_{\text{rel}, \mathbb{E}} = \left| \frac{\mathbb{E}\left[g\right] - \widetilde{\mathbb{E}\left[g\right]}}{\mathbb{E}\left[g\right]} \right|,  
\end{equation}
where $\mathbb{E}\left[g\right]$ is a reference value and $\widetilde{\mathbb{E}\left[g\right]}$ is an estimate computed by directly post-processing the terms of $\widetilde{g}$.
In all numerical experiments, the reference expected value is computed with a quasi-\gls{mc} integration scheme based on Sobol sequences, where the size of the quasi-\gls{mc} sample is equal to $10^8$.

\subsection{Accuracy versus Costs}
\label{subsec:costs}
In the following numerical experiments, the performance of both polynomial approximations is measured by the respective error-cost relation.  
In most engineering applications, the cost of running simulations for different sets of parameter values typically outweighs any other numerical operation necessary for the construction of the approximation. 
In this context, the cost of a method refers solely to the number of simulations, equivalently, solver calls or model evaluations, that are needed until the sought approximation accuracy is reached.
This notion of cost is also used in this section, despite the fact that computationally inexpensive models are employed.
In particular, in each numerical experiment and for both methods, we compute $19$ polynomial approximations corresponding to increasing simulation budgets, i.e., for $B = 10, 20, \dots, 100, 200, \dots, 1000$, and compare the errors of Equations  \eqref{eq:cverr_rms} and \eqref{eq:relerr_mean} for the same costs.

We should note that regression-based \gls{gpc} methods, such as the \gls{lar}-\gls{gpc} \cite{blatman2011} method employed here for comparison purposes, are significantly more expensive than interpolation-based methods in terms of the neglected costs.
This is due to the fact that the solution of possibly expensive \gls{ls} problems is necessary.
On the other hand, regression-based \gls{gpc} methods allow to pre-compute the \gls{qoi} for different parameter values, i.e., run the required simulations offline before the algorithm for the construction of the \gls{gpc} approximation is used.
Moreover, this offline task is embarrassingly parallelizable.
On the contrary, as can easily be observed from Algorithm~\ref{algo:dimadapt}, dimension-adaptive sparse interpolation is based on sequential computations and therefore does not allow the use of parallelization to that extent.
Despite it not being addressed here, those differences must be taken into account by practitioners, considering the problem and the computational resources at hand.

\subsection{Borehole Model}
\label{subsec:borehole}
We consider the $8$-dimensional parametric function:
\begin{equation}
\label{eq:borehole}
g\left(\mathbf{y}\right) = \frac{2 \pi T_{\text{u}} \left(H_{\text{u}} - H_{\text{l}}\right)}{\ln\left({\frac{r}{r_{\text{w}}}}\right) \left( 1 + \frac{T_{\text{u}}}{T_{\text{l}}} + \frac{2 L T_{\text{u}}}{\ln\left({\frac{r}{r_{\text{w}}}}\right) r_{\text{w}}^2 K_{\text{w}}}\right)},
\end{equation}
which models the water flow through a borehole \cite{surjanovic2019}.
The water flow is given in m$^3$/y and the parameter vector is $\mathbf{y} = \left(r_{\text{w}}, r, T_{\text{u}}, H_{\text{u}}, T_{\text{l}}, H_{\text{l}}, L, K_{\text{w}}\right)$,  {where $r_\text{w}$ and $r$ respectively denote the radius of the borehole and the radius of influence (in m), $T_{\text{u}}$ and $T_{\text{l}}$ the transmissivity of the upper and lower aquifer (in m$^2$/y), $H_{\text{u}}$ and $H_{\text{l}}$ the potentiometric head of the upper and lower aquifer (in m), $L$ the length of the borehole (in m), and $K_{\text{w}}$ the hydraulic conductivity of the borehole (in m/y).}

%\begin{table}[b!]
%	\caption{Input parameters of the borehole model.}
%	\centering\begin{tabular}{|c|c|c|}
%		\hline
%		\textbf{Parameter} & \textbf{Symbol} & \textbf{Units} \\
%		\hline
%		radius of borehole & $r_\text{w}$ & m  \\
%		\hline 
%		radius of influence & $r$ & m  \\
%		\hline
%		transmissivity of upper aquifer & $T_{\text{u}}$ & m$^2$/yr  \\
%		\hline
%		potentiometric head of upper aquifer & $H_{\text{u}}$ & m \\
%		\hline
%		transmissivity of lower aquifer & $T_{\text{l}}$ & m$^2$/yr  \\
%		\hline
%		potentiometric head of lower aquifer & $H_{\text{l}}$ & m  \\
%		\hline
%		length of borehole & $L$ & m \\
%		\hline
%		hydraulic conductivity of borehole & $K_{\text{w}}$ & m/yr \\
%		\hline
%	\end{tabular}
%	\label{tab:borehole_params}
%\end{table}

 {Regarding the parameter distributions, we} follow the setting given in \cite{surjanovic2019}, where specific value ranges and probability distributions are given for each parameter and use this information to recast all parameters as to follow truncated normal distributions.
In this numerical experiment, the input parameters are modeled as follows.
\begin{itemize}
	\item The parameter $r_{\text{w}}$ originally follows the normal distribution $\mathcal{N}\left(\mu=0.1, \sigma=0.0161812^2\right)$.
	The distribution is now truncated to the range $\left[0.05, 0.15\right]$;
	\item The parameter $r$ originally follows the log-normal distribution $\mathcal{LN}\left(\mu_{\text{LN}}=7.71, \sigma_{\text{LN}}=1.0056\right)$.
	Therefore, the parameter $r$ has a mean value equal to $\exp\left(\mu_{\text{LN}} + \frac{\sigma_{\text{LN}}^2}{2}\right)$ and a variance equal to $\left(\exp\left(\sigma_{\text{LN}}^2\right) - 1\right) \exp\left(2\mu_{\text{LN}} + \sigma_{\text{LN}}^2\right)$.
	The corresponding truncated normal distribution is defined by these mean and variance values, as well as by the truncation range $\left[100, 50000\right]$l
	\item The remaining parameters are originally uniformly distributed,  {such that $T_{\text{u}} \in \left[63070, 115600\right]$, $H_{\text{u}} \in \left[990, 1110\right]$, $T_{\text{l}} \in \left[63.1, 116\right]$, $H_{\text{l}} \in \left[700, 820\right]$, $L \in \left[1120, 1680\right]$, and $K_{\text{w}} \in \left[9855, 12045\right]$.}
	Assuming a uniform distribution with support in $\left[a, b\right]$, the corresponding truncated normal distribution is given as $\mathcal{TN}\left(\mu=\frac{a+b}{2}, \sigma^2=\frac{(b-a)^2}{12}, l=a, u=b\right)$, i.e., the mean value and variance of the normal distribution correspond to those of the original uniform distribution, while the truncation limits coincide with the uniform distribution's support boundaries. 
\end{itemize}

\begin{figure}[H]
	\begin{tabular}[b]{c}
		\begin{tikzpicture}
		\begin{semilogyaxis}[width=0.45\textwidth, xlabel=Function calls, ylabel=$\epsilon_{\text{RMS}}$, legend pos=north east, grid=both]
		\addplot[mark=None, black, thick] table[x index=0, y index=3]{plot_data/borehole_leja.txt};
		\addplot[mark=None, black, dashed, thick] table[x index=0, y index=4, col sep=comma]{plot_data/borehole_lar.txt};
		\legend{Leja interpolation, LAR-gPC}
		\end{semilogyaxis}
		\end{tikzpicture} \\
		{(\textbf a)}
	\end{tabular}
	\hfill
	\begin{tabular}[b]{c}
		\begin{tikzpicture}
		\begin{semilogyaxis}[width=0.45\textwidth, xlabel=Function calls, ylabel=$\epsilon_{\text{rel}, \mathbb{E}}$, legend pos=north east, grid=both]
		\addplot[mark=None, black, thick] table[x index=0, y index=4]{plot_data/borehole_leja.txt};
		\addplot[mark=None, black, dashed, thick] table[x index=0, y index=5, col sep=comma]{plot_data/borehole_lar.txt};
		\legend{Leja interpolation, LAR-gPC}
		\end{semilogyaxis}
		\end{tikzpicture} \\
		{(\textbf b)}
	\end{tabular}
	\caption{Cost-error relation for the approximations of the borehole model. The approximations are constructed with a dimension-adaptive weighted Leja interpolation algorithm \cite{loukrezis2019} and with a degree-adaptive \gls{lar}-\gls{gpc} (least angle regression-generalized polynomial chaos) algorithm \cite{blatman2011, marelli2015}. The size of the validation sample is $Q=10^5$. The reference expected value is computed with quasi-\gls{mc} (Monte Carlo) integration based on a Sobol sample size equal to $10^8$: (\textbf a) \gls{rms} (root-mean-square) validation error, (\textbf b) Expected value relative error.}
	\label{fig:borehole}
\end{figure}
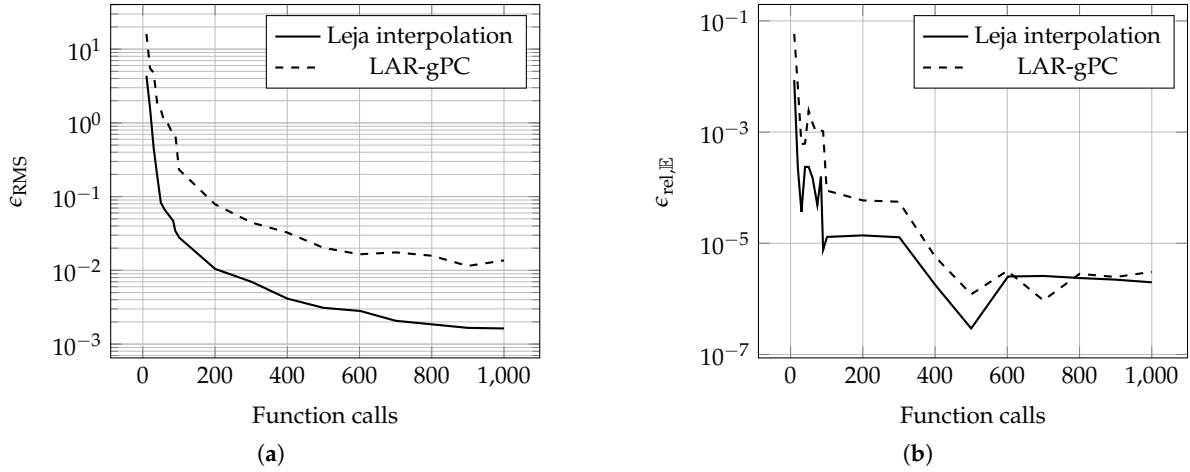

The borehole model of Equation  \eqref{eq:borehole} is approximated by both the dimension-adaptive weighted Leja interpolation and the degree-adaptive \gls{lar}-\gls{gpc} expansion.
Figure~\ref{fig:borehole} shows the error-cost relation for both approximations, with respect to the errors in Equations \eqref{eq:cverr_rms} and \eqref{eq:relerr_mean}.
Regarding the $\epsilon_{\text{RMS}}$ error, Leja interpolation outperformed the \gls{lar}-\gls{gpc} method by approximately one order of magnitude over the whole cost-range.
Leja interpolation was also advantageous with respect to the $\epsilon_{\text{rel}, \mathbb{E}}$ error, again always outperforming the \gls{lar}-\gls{gpc} method and stagnating much sooner.
We note that this error stagnation was observed due to the fact that both methods reached the accuracy of the quasi-\gls{mc} estimate.
Therefore, for the case of the borehole model, the weighted Leja interpolation approach was not only found to perform well for the considered truncated normal input \glspl{pdf}, but was also superior to the \gls{lar}-\gls{gpc} method for both considered error metrics.

\subsection{Steel Column Limit State Function}
\label{subsec:steel}
We consider the $10$-dimensional parametric limit state function:
\begin{equation}
\label{eq:steel_column} 
g\left(\mathbf{y}\right) = F_\text{s} - P_{\text{t}} \left( \frac{1}{2 B D} + \frac{F_0 E_\text{b}}{B D H \left(E_\text{b} - P_{\text{t}}\right)}\right),
\end{equation}
which relates the reliability of a steel column to its cost \cite{kuschel1997, surjanovic2019}.
 {On the right hand side of Equation~\eqref{eq:steel_column}, $P_{\text{t}} = P_\text{d} + P_1 + P_2$ is the total load and $E_{\text{b}} = \frac{\pi^2 E B D H^2}{2 L^2}$ is known as the Euler buckling load.
Thereby, the parameter vector $\mathbf{y}$ consists of ten parameters, namely the yield stress $F_{\text{s}}$ (in MPa), the dead weight load $P_\text{d}$ (in N), the variable loads $P_1$ and $P_2$ (in N), the flange breadth $B$ (in mm), the flange thickness $D$ (in mm), the profile height $H$ (in mm), the initial deflection $F_0$ (in mm), Young's modulus $E$ (in MPa), and the column length $L$ (in mm).}

%\begin{table}[t!]
%	\caption{Input parameters of the steel column function.}
%	\centering\begin{tabular}{|c|c|c|}
%		\hline
%		\textbf{Parameter} & \textbf{Symbol} & \textbf{Units} \\
%		\hline
%		yield stress & $F_\text{s}$ & MPa \\
%		\hline 
%		dead weight load & $P_\text{d}$ & N\\
%		\hline
%		variable load 1 & $P_1$ & N \\
%		\hline
%		variable load 2 & $P_2$ & N \\
%		\hline
%		flange breadth & $B$ & mm \\
%		\hline
%		flange thickness & $D$ & mm \\
%		\hline
%		profile height & $H$ & mm \\
%		\hline  
%		initial deflection & $F_0$ & mm \\
%		\hline
%		Young's modulus & $E$ & MPa \\
%		\hline
%		column length & $L$ & mm \\
%		\hline  
%	\end{tabular}
%	\label{tab:steel_column_params}
%\end{table}

{In the original setting given in \cite{surjanovic2019, kuschel1997}, $P_1$ and $P_2$ are modeled with Gumbel distributions, $E$ follows a Weibull distribution, $F_{\text{s}}$, $B$, $D$, and $H$ are log-normally distributed, and $P_\text{d}$ and $F_0$ are normally distributed.
In this numerical experiment, we replace all log-normal and normal distributions by truncated normal distributions and use the Gumbel distribution to model the remaining three parameters.}
The corresponding distributions are given in Table~\ref{tab:steel_column_dists}.  
With the exception of the parameter $L$, for which the truncation limits are set to $\mu \pm 4 \sigma$, all other truncated normal distributions have support in $\left[\mu - 3 \sigma, \mu + 3 \sigma\right]$.
For all parameters that follow a truncated normal distribution, $\mu$ and $\sigma^2$ coincide with the original mean and variance values given in \cite{surjanovic2019, kuschel1997}.
Accordingly, the  {location and scale parameters} of the Gumbel distributions are chosen such that the original mean values and variances given in \cite{surjanovic2019, kuschel1997} are preserved.

\begin{table}[H]
	\caption{Parameter distributions of the steel column function.}
	\centering\begin{tabular}{cc}
		\toprule
		\textbf{Parameter} & \textbf{Distribution} \\
		\midrule
		$F_\text{s}$ & $\mathcal{TN}\left(\mu=400, \sigma^2=35^2, l=295,u=505\right)$  \\
		%\hline %According to reviewer's comment: In Table 1 "F" should be "Fs" that is the "s" subscript is missing, so we added. Please confirm.
		$P_\text{d}$ & $\mathcal{TN}\left(\mu=500000, \sigma^2=50000^2, l=350000,u=650000\right)$  \\
		%\hline
		$P_1$ & $\mathcal{G}\left(\ell=559495, \beta=70173\right)$  \\
		%\hline
		$P_2$ & $\mathcal{G}\left(\ell=559495, \beta=70173\right)$ \\
		%\hline
		$B$ & $\mathcal{TN}\left(\mu=300, \sigma^2=3^2, l=291, u=309\right)$  \\
		%\hline
		$D$ & $\mathcal{TN}\left(\mu=20, \sigma^2=2^2, l=14, u=26\right)$  \\
		%\hline
		$H$ & $\mathcal{TN}\left(\mu=300, \sigma^2=5^2, l=285, u=315\right)$ \\
		%\hline
		$F_0$ & $\mathcal{TN}\left(\mu=30, \sigma^2=10^2, l=0, u=60\right)$ \\
		%\hline
		$E$ & $\mathcal{G}\left(\ell=208110, \beta=3275\right)$ \\
		%\hline
		$L$ & $\mathcal{TN}\left(\mu=7500, \sigma^2=7.5^2, l=7470, u=7530\right)$ \\
		\bottomrule
	\end{tabular}
	\label{tab:steel_column_dists}
\end{table}

\begin{figure}[H]
	\begin{tabular}[b]{c}
		\begin{tikzpicture}
		\begin{semilogyaxis}[width=0.45\textwidth, xlabel=Function calls, ylabel=$\epsilon_{\text{RMS}}$, legend pos=north east, grid=both]
		\addplot[mark=None, black, thick] table[x index=0, y index=3]{plot_data/steel_column_leja.txt};
		\addplot[mark=None, black, dashed, thick] table[x index=0, y index=4, col sep=comma]{plot_data/steel_column_lar.txt};
		\legend{Leja interpolation, LAR-gPC}
		\end{semilogyaxis}
		\end{tikzpicture} \\
		{(\textbf a)}
	\end{tabular}
	\hfill
	\begin{tabular}[b]{c}
		\begin{tikzpicture}
		\begin{semilogyaxis}[width=0.45\textwidth, xlabel=Function calls, ylabel=$\epsilon_{\text{rel}, \mathbb{E}}$, legend pos=north east, grid=both]
		\addplot[mark=None, black, thick] table[x index=0, y index=4]{plot_data/steel_column_leja.txt};
		\addplot[mark=None, black, dashed, thick] table[x index=0, y index=5, col sep=comma]{plot_data/steel_column_lar.txt};
		\legend{Leja interpolation, LAR-gPC}
		\end{semilogyaxis}
		\end{tikzpicture} \\
		{(\textbf b)}
	\end{tabular}
	\caption{Cost-error relation for the approximations of the steel column function. The approximations are constructed with a dimension-adaptive weighted Leja interpolation algorithm \cite{loukrezis2019} and with a degree-adaptive \gls{lar}-\gls{gpc} algorithm \cite{blatman2011, marelli2015}. The size of the validation sample is $Q=10^5$. The reference expected value is computed with quasi-\gls{mc} integration based on a Sobol sample size equal to $10^8$: (\textbf a) \gls{rms} validation error, (\textbf b) Expected value relative error.}
	\label{fig:steel_column}
\end{figure}
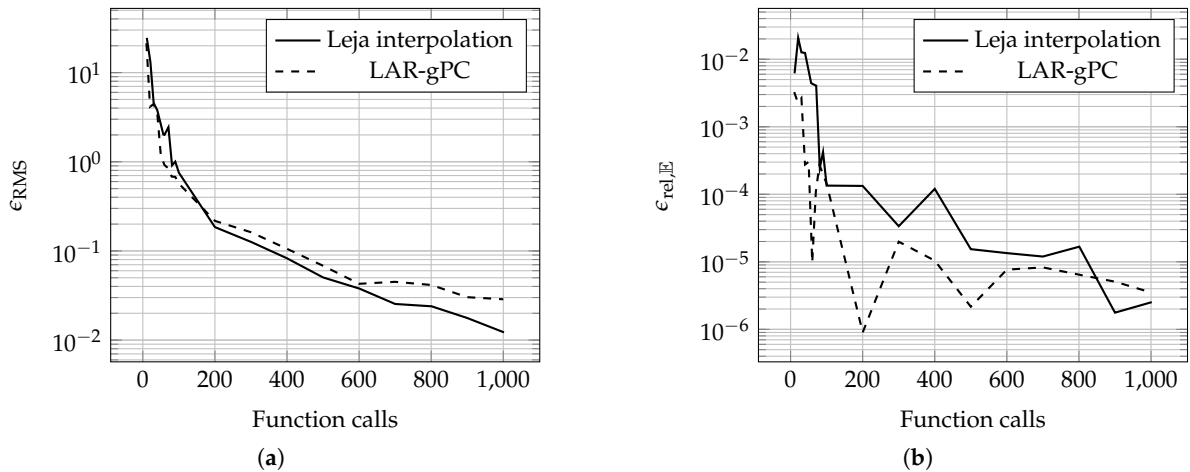

We approximate the steel column function Equation \eqref{eq:steel_column} using the dimension-adaptive weighted Leja interpolation and the degree-adaptive, \gls{lar}-based \gls{gpc} method.
For both approximations and for both error metrics in  Equations \eqref{eq:cverr_rms} and \eqref{eq:relerr_mean}, the error-cost relation is presented in Figure~\ref{fig:steel_column}. 
Once more, the Leja interpolation method performs very well for the given input \glspl{pdf}.
However, contrary to the results of Section~\ref{subsec:borehole}, no advantage is observed over the \gls{lar}-\gls{gpc} method.
With respect to the $\epsilon_{\text{RMS}}$ error, the performance of both methods is comparable.
This can be probably attributed to a reduced regularity of the \gls{qoi} defined in Equation \eqref{eq:steel_column}, in comparison to the \gls{qoi} given in Equation \eqref{eq:borehole}.
The Leja interpolation method has a slight edge for costs greater than $200$ simulations, however, the difference is minor.
With the exception of costs greater than $800$ simulations, the \gls{lar}-\gls{gpc} method is superior to the Leja interpolation in terms of the $\epsilon_{\text{rel}, \mathbb{E}}$ error.
Overall, the weighted Leja interpolation is again able to yield accurate approximations and statistics and it may be regarded as competitive to the \gls{lar}-\gls{gpc} method.

\subsection{Meromorphic Function}
\label{subsec:meromorphic}
We consider the $16$-dimensional meromorphic function:
\begin{equation}
\label{eq:meromorphic}
g\left(\mathbf{y}\right) = \frac{1}{1 + \mathbf{w} \cdot \mathbf{y}},
\end{equation}
where $\mathbf{w}$ is a vector of positive weights.
The weight vector in Equation  \eqref{eq:meromorphic} is given by $\mathbf{w} = \frac{\hat{\mathbf{w}}}{2 \|\hat{\mathbf{w}}\|_1}$, where $\hat{\mathbf{w}} = \left(1, 0.5, 0.1, 0.05, 0.01, \dots, 5\cdot 10^{-8}\right)$ \cite{migliorati2013a}.
The input vector $\mathbf{y}$ takes values in the image space of the random vector $\mathbf{Y}=\left(Y_1, Y_2, \dots, Y_{16}\right)$, where the single \glspl{rv} are given as:
\begin{align*}
\label{eq:meromorphic_dists}
Y_n \sim  
\begin{cases}
\mathcal{TN}\left(\mu=0, \sigma^2=1, l=0, u=3\right), &n=1,3,\dots,15, \\
\mathcal{TN}\left(\mu=0, \sigma^2=1, l=-3, u=0\right), &n=2,4,\dots,16.
\end{cases}
\end{align*}
While a normal distribution truncated exactly at its {(untruncated)} mean value is not very likely to be encountered in practical applications, the selected truncation limits result in very atypical \glspl{pdf} which fit very well within the concept of arbitrary input \glspl{pdf}.

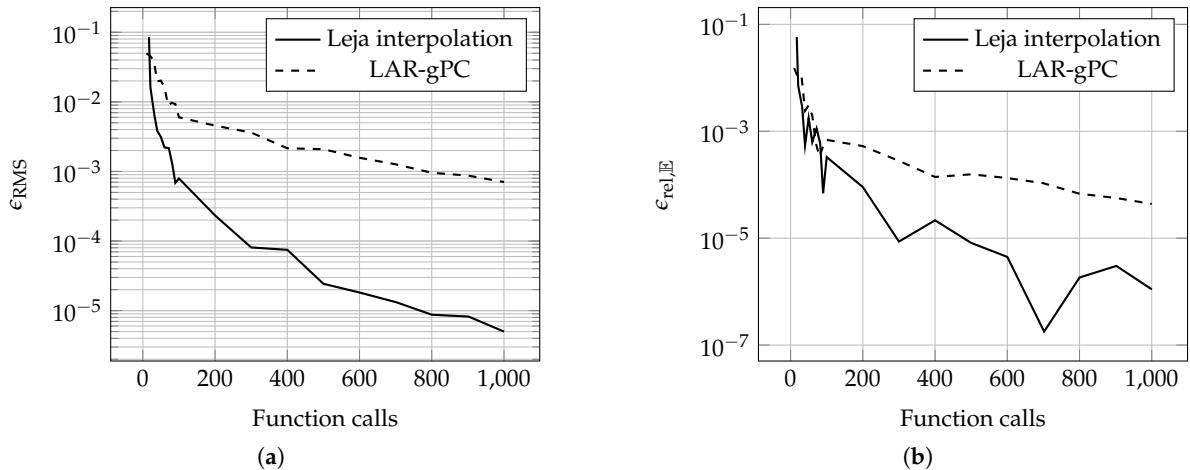
\begin{figure}[H]
	\begin{tabular}[b]{c}
		\begin{tikzpicture}
		\begin{semilogyaxis}[width=0.45\textwidth, xlabel=Function calls, ylabel=$\epsilon_{\text{RMS}}$, legend pos=north east, grid=both]
		\addplot[mark=None, black, thick] table[x index=0, y index=3]{plot_data/meromorphic_leja.txt};
		\addplot[mark=None, black, dashed, thick] table[x index=0, y index=4, col sep=comma]{plot_data/meromorphic_lar.txt};
		\legend{Leja interpolation, LAR-gPC}
		\end{semilogyaxis}
		\end{tikzpicture} \\
		{(\textbf a)}
	\end{tabular}
	\hfill
	\begin{tabular}[b]{c}
		\begin{tikzpicture}
		\begin{semilogyaxis}[width=0.45\textwidth, xlabel=Function calls, ylabel=$\epsilon_{\text{rel}, \mathbb{E}}$, legend pos=north east, grid=both]
		\addplot[mark=None, black, thick] table[x index=0, y index=4]{plot_data/meromorphic_leja.txt};
		\addplot[mark=None, black, dashed, thick] table[x index=0, y index=5, col sep=comma]{plot_data/meromorphic_lar.txt};
		\legend{Leja interpolation, LAR-gPC}
		\end{semilogyaxis}
		\end{tikzpicture} \\
		{(\textbf b)}
	\end{tabular}
	\caption{Cost-error relation for the approximations of the meromorphic function. The approximations are constructed with a dimension-adaptive weighted Leja interpolation algorithm \cite{loukrezis2019} and with a degree-adaptive \gls{lar}-\gls{gpc} algorithm \cite{blatman2011, marelli2015}. The size of the validation sample is $Q=10^5$. The reference expected value is computed with quasi-\gls{mc} integration based on a Sobol sample size equal to $10^8$: (\textbf a) \gls{rms} validation error, (\textbf b) Expected value relative error.}
	\label{fig:meromorphic}
\end{figure}

Figure~\ref{fig:meromorphic} shows the error-cost relations corresponding to the error metrics of Equations \eqref{eq:cverr_rms} and \eqref{eq:relerr_mean} for both approximations of the meromorphic function in Equation \eqref{eq:meromorphic}, i.e., using the dimension-adaptive weighted Leja interpolation algorithm and the degree-adaptive \gls{lar}-\gls{gpc} algorithm.
The results resemble those of Section~\ref{subsec:borehole}, i.e., the Leja interpolation method had a clear advantage over the \gls{lar}-\gls{gpc} approximation with respect to both error metrics.
Regarding the $\epsilon_{\text{RMS}}$ error, the difference between the two approximations was typically greater than one order of magnitude, while it increased for an increasing computational budget.
Regarding the $\epsilon_{\text{rel}, \mathbb{E}}$ error, the two methods showed a comparable performance for costs up to $100$ simulations, however, the Leja interpolation was again clearly superior as the simulation budget increased.
Hence, in this numerical experiment, the weighted Leja interpolation was able to accurately approximate a model with very atypical input \glspl{pdf}, and provide accurate estimations of statistics.
Moreover, it was found to have a clear edge over the well-established \gls{lar}-\gls{gpc} method.

\subsection{Dielectric Inset Waveguide}
\label{subsec:waveguide}
\begin{figure}[H]
	\centering
	\begin{tikzpicture}
	\begin{axis}[
	axis lines = center,
	width=10cm,height=8cm,
	xtick={0},ytick={0},
	%  minor tick={-12,-11,...,12},
	xmin=-2,xmax=20,ymin=-4,ymax=34,
	xlabel={$z$},ylabel={$x$}
	]
	\addplot [no marks, dashed, fill=gray!20] coordinates {(5,0) (12,0) (12,30) (5,30)};
	\addplot [no marks, dashed, fill=white] coordinates {(0,0) (5,0) (5,30) (0,30)};
	\addplot [no marks, fill=white] coordinates {(12,0) (17,0) (17,30) (12,30)};
	\addplot [no marks, dashed] coordinates {(12,0) (12,30)};
	\addplot [no marks] coordinates {(0,0) (0,30)};
	\addplot [no marks] coordinates {(0,0) (20,0)};
	% label points
	\node [left ] at (axis cs:4,15) {$\mu_0,\varepsilon_0$};
	\node [left ] at (axis cs:9.5,15) {$\mu,\varepsilon$};
	\node [left ] at (axis cs:16,15) {$\mu_0,\varepsilon_0$};
	\node [left ] at (axis cs:0,30) {$w$};
	\node [below ] at (axis cs:5,0) {$d$};
	\node [below ] at (axis cs:12,0) {$d+\ell$};
	\node [below ] at (axis cs:17,0) {$2d+\ell$};
	% ports
	\addplot [no marks, solid, black, line width=2pt] coordinates {(0,0) (0,30)};
	\node [left ] at (axis cs:0,17) {$\Gamma_{\text{in}}$};
	\addplot [no marks, solid, black, line width=2pt] coordinates {(17,0) (17,30)};
	\node [right ] at (axis cs:17,17) {$\Gamma_{\text{out}}$};
	% PEC
	\addplot [no marks, solid, black!50, line width=2pt] coordinates {(0,30) (17,30)};
	\addplot [no marks, solid, black!50, line width=2pt] coordinates {(0,0) (17,0)};
	\node [above ] at (axis cs:9,30) {$\Gamma_{\text{PEC}}$};
	\end{axis}
	\end{tikzpicture}
	\caption{View of the waveguide in the $xz$ plane. The input and output ports are shown in black, the \gls{pec} (perfect electric conductors) walls in dark gray, the vacuum-filled area in white, and the dielectric material in light gray.}
	\label{fig:xzview}
\end{figure}
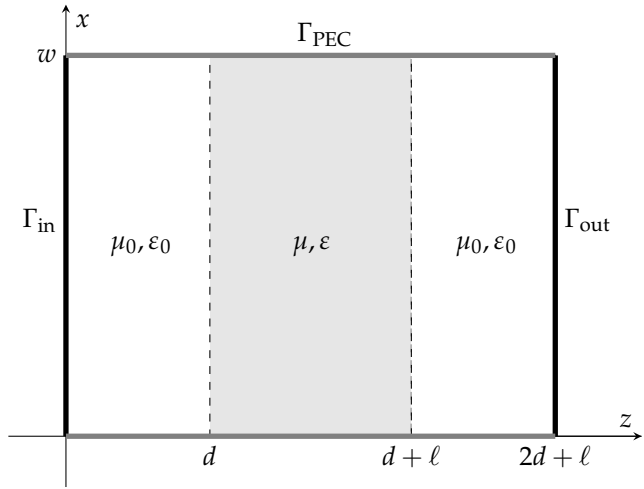

We consider a rectangular waveguide with a dispersive dielectric inset, similar to the one investigated in \cite{loukrezis2019c}.
A 2D illustration of the waveguide is depicted in Figure~\ref{fig:xzview}, where $w$ denotes its width, $\ell$ is the length of the dielectric material, and $d$ is a vacuum offset. 
The waveguide is further extended in the $y$ direction by its height $h$, which is not shown in Figure~\ref{fig:xzview}.
The dielectric material has permittivity $\varepsilon = \varepsilon_0 \varepsilon_{\mathrm{r}}$ and permeability $\mu = \mu_0 \mu_{\mathrm{r}}$, where the subscripts ``$0$'' and ``$\mathrm{r}$'' refer to the absolute value of the material property in vacuum and to its relative value for the given dielectric material, respectively.
The relative material values are given by Debye relaxation models of second order \cite{xu2010}, such that:
\begin{align*}
\varepsilon_{\mathrm{r}} &= \varepsilon_{\infty} + \frac{\varepsilon_{\mathrm{s},1} - \varepsilon_{\infty}}{1 + \left(\imath \omega \tau_{\varepsilon,1}\right)} + \frac{\varepsilon_{\mathrm{s},2} - \varepsilon_{\infty}}{1 + \left(\imath \omega \tau_{\varepsilon,2}\right)}, \\  
\mu_{\mathrm{r}} &= \mu_{\infty} + \frac{\mu_{\mathrm{s},1} - \mu_{\infty}}{1 + \left(\imath \omega \tau_{\mu,1}\right)} + \frac{\mu_{\mathrm{s},2} - \mu_{\infty}}{1 + \left(\imath \omega \tau_{\mu,2}\right)},
\end{align*}
where $\tau_{\varepsilon/\mu,1/2}$ are relaxation time constants, the subscript ``$\infty$'' refers to a very high frequency value of the relative material property, the subscript ``s'' to a static value of the relative material property, and $\imath$ denotes the imaginary unit.
The waveguide's input and output boundaries are denoted with $\Gamma_{\text{in}}$ and $\Gamma_{\text{out}}$, respectively.
The remaining waveguide walls are assumed to be \glspl{pec} and the corresponding boundary is denoted with $\Gamma_{\text{PEC}}$.

We assume that the input port boundary $\Gamma_\text{\text{in}}$ is excited by an incoming plane wave. 
We further assume that the excitation coincides with the fundamental \gls{te} mode and that all other propagation modes can be neglected, e.g., are quickly attenuated.
Under these assumptions, the electric field $\mathbf{E}$ can be computed by solving the boundary value problem: 
\begin{subequations}
	\label{eq:maxwell_source}
	\begin{align}
	\nabla \times \left( \mu^{-1}\nabla \times \mathbf{E}\right) - \omega^2 \varepsilon \mathbf{E} &= 0, \quad &&\text{~in~} D,\\
	\mathbf{n}_{\Gamma_{\text{PEC}}} \times \mathbf{E} &= 0, \quad  &&\text{~on~} \Gamma_{\text{PEC}}, \\
	\mathbf{n}_{\Gamma_{\text{in}}} \times \left(\nabla \times \mathbf{E}\right) + \imath k_z \mathbf{n}_{\Gamma_{\text{in}}} \times \left(\mathbf{n}_{\Gamma_{\text{in}}} \times \mathbf{E}\right) &= \mathbf{U}^{\text{i}}, \quad  &&\text{~on~} \Gamma_{\text{in}}, \\
	\mathbf{n}_{\Gamma_{\text{out}}} \times \left(\nabla \times \mathbf{E}\right) + \imath k_z \mathbf{n}_{\Gamma_{\text{out}}} \times \left(\mathbf{n}_{\Gamma_{\text{out}}} \times \mathbf{E}\right) &= 0, \quad  &&\text{~on~} \Gamma_{\text{out}},
	\end{align}
	\label{eq:Maxwell_boundary_problem}
\end{subequations}
where $D$ is the computational domain, $\omega = 2 \pi f$ the angular frequency, $f$ the frequency, $\mathbf{U}^{\text{i}}$ the incoming plane wave, $\mathbf{n}_{\Gamma_{\text{in}}/\Gamma_{\text{out}}/\Gamma_{\text{PEC}}}$ are outwards-pointing normal vectors, and $\mathbf{k} = \left(0,0,k_z\right)$ the wavevector.
The \gls{qoi} is the reflection coefficient at the input port $\Gamma_\text{\text{in}}$,  $r = \left| \frac{\mathbf{E}_{\Gamma_{\text{in}}^-}}{\mathbf{E}_{\Gamma_{\text{in}}^+}}\right| \in \left[0,1\right]$. 
Usually, Equation \eqref{eq:maxwell_source} is solved numerically, e.g., using the \gls{fem}. 
For this simple model, a  semi-analytical solution exists for the reflection coefficient $r$. 
Therefore, errors due to spatial discretization can be neglected.

Equation~\eqref{eq:Maxwell_boundary_problem} is transformed into a stochastic problem by recasting all geometrical and Debye material model parameters, as well as the frequency, to \glspl{rv}.
In the nominal configuration, the parameter values are $\bar{f}=\SI{6}{\giga\hertz}$, $\bar{w}=\SI{30}{\milli\meter}$, $\bar{h}=\SI{3}{\milli\meter}$, $\bar{\ell}=\SI{7}{\milli\meter}$, $\bar{d}=\SI{1}{\milli\meter}$, $\bar{\varepsilon}_{\text{s},1}=2$, $\bar{\varepsilon}_{\text{s},2}=2.2$, $\bar{\varepsilon}_{\infty}=1$,  $\bar{\mu}_{\text{s},1}=2$, $\bar{\mu}_{\text{s},2}=3$, $\bar{\mu}_{\infty}=1$, $\bar{\tau}_{\varepsilon,1} = 1$, $\bar{\tau}_{\varepsilon,2} = 1.1$, $\bar{\tau}_{\mu,1} = 1$, and $\bar{\tau}_{\mu,2} = 2$.
Each parameter is now assumed to follow a truncated normal distribution, such that $Y_n \sim \mathcal{TN}\left(\bar{y}_n, \left(0.01\,\bar{y}_n\right)^2, 0.95 \, \bar{y}_n, 1.05 \, \bar{y}_n\right)$, $n=1,\dots,15$, i.e., the standard deviation of the Gaussian distribution is equal to $1\%$ of the nominal parameter value and a maximum deviation of $\pm 5\%$ from the nominal parameter value is imposed.

\begin{figure}[H]
	\begin{tabular}[b]{c}
		\begin{tikzpicture}
		\begin{semilogyaxis}[width=0.45\textwidth, xlabel=Function calls, ylabel=$\epsilon_{\text{RMS}}$, legend pos=north east, grid=both]
		\addplot[mark=None, black, thick] table[x index=0, y index=2]{plot_data/cv_results_wg_debye2_dali.txt};
		\addplot[mark=None, black, dashed, thick] table[x index=0, y index=3, col sep=comma]{plot_data/wg_debye2_lar.txt};
		\legend{Leja interpolation, LAR-gPC}
		\end{semilogyaxis}
		\end{tikzpicture} \\
		{(\textbf a)}
	\end{tabular}
	\hfill
	\begin{tabular}[b]{c}
		\begin{tikzpicture}
		\begin{semilogyaxis}[width=0.45\textwidth, xlabel=Function calls, ylabel=$\epsilon_{\text{rel}, \mathbb{E}}$, legend pos=north east, grid=both]
		\addplot[mark=None, black, thick] table[x index=0, y index=1]{plot_data/moments_relerr_wg_debye2_dali.txt};
		\addplot[mark=None, black, dashed, thick] table[x index=0, y index=4, col sep=comma]{plot_data/wg_debye2_lar.txt};
		\legend{Leja interpolation, LAR-gPC}
		\end{semilogyaxis}
		\end{tikzpicture} \\
		{(\textbf b)}
	\end{tabular}
	\caption{Cost-error relation for the approximations of the dielectric inset waveguide model. The approximations are constructed with a dimension-adaptive weighted Leja interpolation algorithm \cite{loukrezis2019} and with a degree-adaptive \gls{lar}-\gls{gpc} algorithm \cite{blatman2011, marelli2015}. The size of the validation sample is $Q=10^5$. The reference expected value is computed with quasi-\gls{mc} integration based on a Sobol sample size equal to $10^8$: (\textbf a) \gls{rms} validation error, (\textbf b) Expected value relative error.}
	\label{fig:wg_debye2}
\end{figure}
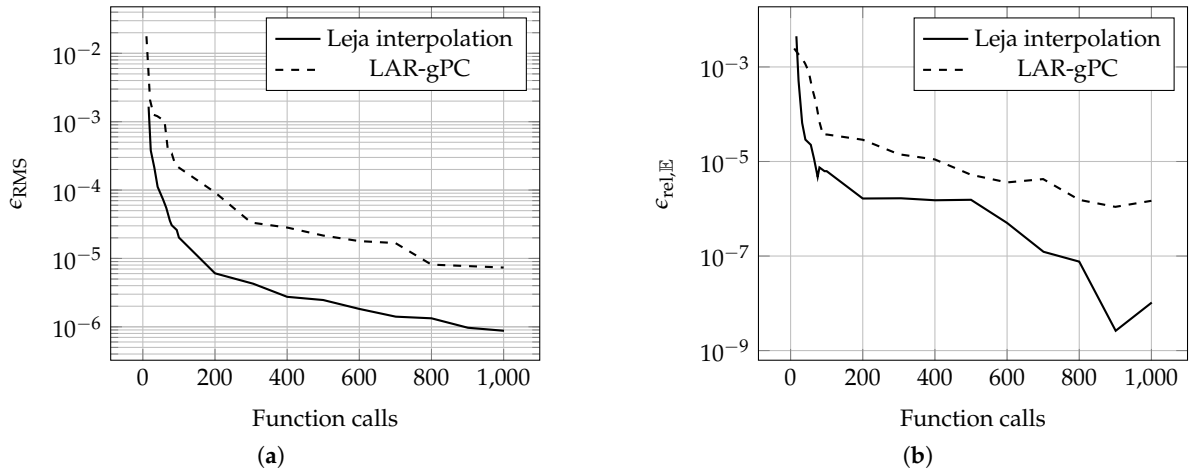

The comparison between the dimension-adaptive Leja interpolation and the degree-adaptive \gls{lar}-\gls{gpc} for both error metrics in Equations \eqref{eq:cverr_rms} and \eqref{eq:relerr_mean} is presented in Figure~\ref{fig:wg_debye2}. 
Similarly to the results obtained in Sections~\ref{subsec:borehole} and \ref{subsec:meromorphic}, the Leja interpolation method was significantly superior to the \gls{lar}-\gls{gpc} approximation.
For both error metrics, the difference was typically one order of magnitude,  with an increasing tendency for larger simulation budgets.

\section{Discussion}
\label{sec:concl}
In this work, we investigated the use of a sparse interpolation method based on weighted Leja node sequences in order to address the problems of approximation and uncertainty quantification for {models with input random variables} characterized by \glspl{pdf} of arbitrary shapes.
For that purpose, we applied a weighted Leja-based, dimension-adaptive interpolation algorithm to four models featuring $8$, $10$, $16$, and $15$ parameters, respectively.
Truncated normal and extreme value input distributions have been used to model the random input parameters. 
The suggested approach has also been compared to a well-established adaptive \gls{gpc} method, which uses numerically constructed orthogonal polynomials to address arbitrary input \glspl{pdf}.

The numerical results presented in Section~\ref{sec:num_exp} indicated that the weighted Leja interpolation performed very well in all four numerical experiments, both in terms of approximation as well as of statistics estimation accuracy.
Moreover, the suggested approach showed either superior or equivalent performance to the \gls{gpc} method used for comparison purposes.
The comparison results could be seen as an extension of those reported in \cite{eldred2009a, eldred2009}, where uniform, normal, and log-normal distributions were considered, to the case of arbitrary input \glspl{pdf}.
Overall, we may conclude that the approach presented in this work posed a viable and reliable alternative to arbitrary \gls{gpc} methods \cite{soize2004, wan2006a, wan2006}.

An obvious limitation of this work is the assumption that the input \glspl{rv} could be described by continuous \glspl{pdf}.
Data-driven approaches, such as those presented in \cite{ahlfeld2016, oladyshkin2012}, are able to construct \gls{gpc} approximations based only on moment values computed after a dataset, without any other assumption regarding a \gls{pdf}, which might rely on subjective interpretations of the data and result in the introduction of severe biases.
The extension of interpolation-based \gls{uq} methods to problems where the input \glspl{pdf} are not explicitly given in closed-form would be an interesting path for further research, to be pursued in a later study.

%%%%%%%%%%%%%%%%%%%%%%%%%%%%%%%%%%%%%%%%%%

%%%%%%%%%%%%%%%%%%%%%%%%%%%%%%%%%%%%%%%%%%
\vspace{6pt} 

\reftitle{References}


%\externalbibliography{yes}
%\bibliography{references_aLeja.bib}
%=====================================
% References, variant B: internal bibliography
%=====================================

%%%%%%%%%%%%%%%%%%%%%%%%%%%%%%%%%%%%%%%%%%
%% optional
%\sampleavailability{Samples of the compounds ...... are available from the authors.}

%% for journal Sci
%\reviewreports{\\
%Reviewer 1 comments and authors’ response\\
%Reviewer 2 comments and authors’ response\\
%Reviewer 3 comments and authors’ response
%}

%%%%%%%%%%%%%%%%%%%%%%%%%%%%%%%%%%%%%%%%%%
\end{document}